 \documentclass[final,3p,times]{elsarticle}

\usepackage{amssymb}
\usepackage{amsmath}        
\usepackage{amstext}   
\usepackage{mathrsfs}

\usepackage{graphics}
\usepackage{graphicx}
\usepackage{rotating}
\usepackage{psfrag}         
\usepackage{color}
\usepackage{url}
\usepackage[dvips]{epsfig}
\usepackage{multirow}
\usepackage{algorithmic}
\usepackage{algorithm,algorithmic}
\usepackage[caption=false]{subfig}
\usepackage{bm}

\def\ds{\displaystyle}

\DeclareSymbolFont{msbm}{U}{msb}{m}{n}
\DeclareMathSymbol{\N}{\mathalpha}{msbm}{'116}
\DeclareMathSymbol{\Rset}{\mathalpha}{msbm}{'122}

\def\ds{\displaystyle}

\newcommand{\bfU}{{\bf U}}

\newcommand{\bfx}{\boldsymbol{x}}
\newcommand{\bfu}{\boldsymbol{u}}
\newcommand{\bfv}{\boldsymbol{v}}

\begin{document}

\begin{frontmatter}
	
\title{Highly efficient NURBS-based isogeometric analysis for coupled nonlinear diffusion-reaction equations with and without advection}
\author[label1]{Ilham Asmouh\corref{corres}}\cortext[corres]{Corresponding author} 
\ead{ilham.asmouh@uibk.ac.at}
\author[label1]{Alexander Ostermann} 
\ead{alexander.ostermann@uibk.ac.at}

\address[label1]{Institut f\"ur Mathematik, Universit\"at Innsbruck, Technikerstra{\ss}e 13,  6020 Innsbruck, Austria}

\begin{abstract}
Nonlinear diffusion-reaction systems model a multitude of physical phenomena. A common situation is biological development modeling where such systems have been widely used to study spatiotemporal phenomena in cell biology. Systems of coupled diffusion-reaction equations are usually subject to some complicated features directly related to their multiphysics nature. Moreover, the presence of advection is source of numerical instabilities, in general, and adds another challenge to these systems. In this study, we propose a NURBS-based isogeometric analysis (IgA) combined with a second-order Strang operator splitting to deal with the multiphysics nature of the problem. The advection part is treated in a semi-Lagrangian framework and the resulting diffusion-reaction equations are then solved using an efficient time-stepping algorithm based on operator splitting. The accuracy of the method is studied by means of a advection-diffusion-reaction system with analytical solution. To further examine the performance of the new method on complex geometries, the well-known Schnakenberg--Turing problem is considered with and without advection. Finally, a Gray--Scott system on a circular domain is also presented. The results obtained demonstrate the efficiency of our new algorithm to accurately reproduce the solution in the presence of complex patterns on complex geometries. Moreover, the new method clarifies the effect of geometry on Turing patterns.
\end{abstract}

\end{frontmatter}
\section{Introduction}
In this paper, given a bounded domain $\Omega \subset \Rset^d$, $d \in \{1,2,3\}$, with a  sufficiently smooth boundary $\Gamma$ and a given time interval $[0,T]$, we are interested in solving the unsteady nonlinear coupled advection-diffusion-reaction equations
\begin{equation}\label{model}
    \begin{aligned}
       \frac{\partial u}{\partial t} + {\bf a}(\bfx,t)\cdot\nabla{u} &=  \nabla \cdot (d_{1}\nabla u) + f(u,v),\quad (\bfx, t) \in  \Omega \times [0,T], \\
       \frac{\partial v}{\partial t} + {\bf a}(\bfx,t)\cdot\nabla{v} &=  \nabla \cdot (d_{2}\nabla v) + g(u,v),\quad (\bfx, t) \in  \Omega \times [0,T], \\
       u({\bfx},0)  &= u_{0}({\bfx}), \quad \bfx \in  \Omega, \\
       v({\bfx},0)  &= v_{0}({\bfx}), \quad \bfx \in  \Omega,
    \end{aligned}
\end{equation}
with appropriate boundary conditions. In biological development modeling, the unknowns $u$ and $v$ denote the concentrations of the activator and the inhibitor respectively, which are advected by the velocity field ${\bf a}(\bfx,t)$. The quantities $d_{1}$ and $d_{2}$ denote the diffusion coefficients associated to $u$ and $v$. The right-hand sides $f$ and $g$ are the reaction functions which link the activator $u$ and the inhibitor $v$. The model problem \eqref{model} arises in many fields of research such as biological development (cf.~the work of Turing), neuroscience, calcium dynamics, see \cite{sneyd1995intercellular,hodgkin1952quantitative,meron1992pattern} among others. Pattern formation in developmental biology, which concerns the study of animals and plants with their growth and development processes, refers to the generation of complex organizations of cell fates in space and time. Within an embryo, this includes the various steps and processes involved to create an organism. Pattern formation is a direct result of interaction between many components and is expected to be nonlinear. Therefore, the biological patterns are mathematically represented by coupled nonlinear diffusion-reaction equations in general. These models express the temporal behavior of the concentrations of reacting and diffusing chemicals. These problems include the Turing-type models. As examples, we mention the Schnakenberg model \cite{schnakenberg1979simple}, the Gierer--Meinhardt model \cite{gierer1972theory}, the Thomas model, the Gray--Scott model \cite{gray1983autocatalytic} and others presented in \cite{barrio2002size,aragon2002turing}.\\

It is well known that the embryos in developing tissues are usually characterized by the complex shapes of the patterns. As a consequence, a computational tool that supports complex geometries is essential. Another challenge lies in the multiphysical nature of this problems. In fact, the nonlinearity of the reaction part constitutes the most important feature characterizing nonlinear coupled diffusion-reaction problems. Therefore, the development of a powerful tool which solves such problems is still taking attention of many research works. The authors in \cite{wu2018two} proposed a two-level linearized compact alternating direction implicit (ADI) scheme to solve the nonlinear Fisher problem in two space dimensions, where fourth-order accuracy is achieved. To deal with the nonlinearity, the authors used Newton's method. In \cite{toro2009ader}, the authors extended the Godunov and ADER approaches to construct finite volume schemes of arbitrary order of accuracy in space and time to solve nonlinear diffusion-reaction problems. Their approach is essentially based on solving Riemann problems and their study is limited to one space dimension. Mixed finite element methods are applied to the same problem in \cite{chen2011two} where a two-grid algorithm is used to linearize the mixed equations. The authors in \cite{liao2006fourth} presented a fourth-order compact algorithm for nonlinear diffusion-reaction equations with Neumann boundary conditions. Concerning nonlinear systems of diffusion-reaction, several contributions have been aimed at developing numerical tools for the simulation of such coupled equations. For example, the authors in \cite{shakeri2011finite} presented a finite volume spectral element method for solving the Schnakenberg--Turing model in two space dimensions, where a Crank--Nicolson technique is adopted to integrate the time derivative. The numerical results are in good agreement with existing studies although the method is limited to regular geometries. In \cite{dehghan2016numerical}, the authors presented a meshless method based on a moving Kriging element free Galerkin (EFG) approach to solve two and three-dimensional Turing patterns, where the Allen--Cahn, the Gray--Scott, and the Brusselator models are considered. To deal with the time derivative, the authors used an explicit fourth-order Runge--Kutta scheme. The same problems are treated in \cite{hepson2021numerical} using quartic trigonometric B-splines.\\

Most of the previous works are limited to simple regular geometries. However, it is known that pattern formation in such complicated nonlinear models is supported by complex media in general \cite{marshall2011origins}. Therefore, an efficient method for complex geometries is essential. Moreover, the introduction of advection in these models is required, since advection is so fundamental in nature, and it is related to the transport of species due to the presence of velocity fields in general realistic cases. In the current study, the spatial derivatives are treated with the isogeometric analysis (IgA) \cite{hughes2005isogeometric}. It is well known that IgA is designed as one of the most powerful tools for the exact representation of the geometries, no matter how complicated they are, with only a few elements. Since 2006, IgA has been used intensively in different areas, such as plasticity, elasticity and structural vibrations \cite{cottrell2009isogeometric,elguedj2008f,auricchio2007fully}, computational fluid dynamics \cite{bazilevs2010large,calo2008simulation,zhang2007patient}, structural mechanics \cite{cottrell2006isogeometric,kiendl2010bending}, among others. The basis functions in IgA are generated from Non Uniform Rational B-Spline functions (NURBS). These functions possess many properties that make them an attractive choice over traditional Lagrange basis functions. Indeed, the extension to high-space dimensions is straightforward thanks to the tensor product property of NURBS. In addition, the NURBS functions offer high inter-element continuity which is absent in the Lagrange basis functions. This is decisive especially in presence of sharp layers. In order to deal with the multiphysical nature of problem \eqref{model}, we propose a second-order accurate Strang splitting technique in which we separate the nonlinear reactive terms from the diffusive ones and integrate each with an appropriate method. Operator splitting offers the possibility to manage each part of the whole system separately using an appropriate technique, which makes it an excellent choice that contributes in reducing the computational costs. In our study, the advective terms are treated in a semi-Lagrangian framework. The resulting diffusion-reaction equations are then treated separately. Specifically, the diffusion terms are discretized using an isogeometric finite element method, where the semi-discrete equations are integrated using the second-order implicit scheme of Gear also known in the literature as backward differentiation formula (BDF2), while the nonlinear reaction terms are integrated with an explicit Runge--Kutta scheme. By doing so, no linearizing approach is required. Moreover, the semi-Lagrangian approach is proven to produce accurate results allowing large time steps while guaranteeing the stability of the solution during time evolution. As a result, the computational costs are reduced. The accuracy of the new method is firstly tested for a advection-diffusion-reaction system with analytical solution. The Schnakenberg--Turing system, in absence of advection, is treated on different complex geometries and promising results are obtained with the formation of the expected patterns as described in the literature. Next, to analyze the effect of advection on Turing-patterns, we carry out the study with a toroidal velocity field. Finally the Gray--Scott model on a circular domain is presented and compared to existing studies.\\

The rest of this paper is structured as follows. Section \ref{sec2} contains a brief introduction to IgA. This includes the definition of NURBS functions and some important relations. In section \ref{sec3}, the novel isogeometric semi-Lagrangian analysis is detailed. This includes also the description of the $L^{2}$-projection approach used in order to deal with the fact that NURBS are not interpolatory. The Strang operator splitting is presented where the formulation of the new numerical approach is detailed for the diffusive terms and the nonlinear reaction terms. Section \ref{sec4} provides the obtained numerical results for different test problems. This includes, the well-known Schnakenberg--Turing and Gray--Scott models. Finally, concluding remarks are given.
\section{A short summary of NURBS functions}\label{sec2}
In the current work, we use NURBS as an analysis tool to discretize the spatial domain. Classical B-splines are very popular for their stability and simplicity. However, the main drawback of conventional B-spline techniques is their inability to accurately fit curves regardless of their degree. In fact, B-spline functions can represent only segments of a straight line and some curves and surfaces (but not most polynomial curves of degree 2 or higher \cite{grove2011ct,rogers2001introduction}). The representation of a plate with a circular hole, for example, requires only two elements in the case of NURBS \cite{hughes2005isogeometric}. This is particularly useful for complex geometries, as it both accurately represents the desired geometry and minimizes the number of elements required to mesh the domain. This reduces the computational cost. We start with a brief summary of B-spline functions and NURBS-based IgA in one and two space dimensions. \\

Let $\Xi = \left(\xi_{1},\xi_{2},\dots,\xi_{m_{b}+p+1}\right)$ be a given knot vector which consists of an ordered set of non-decreasing parameter values, where $\xi_{i},\ 1\le i\le m_{b}+p+1$ represent the knots which partition the parameter interval into elements. The integer $m_{b}$ denotes the number of basis functions of chosen degree $p$. Following \cite{cottrell2009isogeometric,nguyen2015isogeometric}, e.g., the NURBS functions are then derived from B-splines as
\begin{equation}\label{NURBSbasis}
\mathcal{\widetilde{N}}^{i}_{p}(\xi)=\frac{\mathcal{B}^{p}_{i}(\xi)\omega_{i}}{W(\xi)},
\end{equation}
where $\mathcal{B}^{p}_{i}(.)$ is the B-spline of degree $p$, the numbers $\left\{\omega_{i}\right\}_{i=1}^{m_{b}}$ are the NURBS weights and $W(\xi)=\ds\textstyle\sum_{i=1}^{m_{b}}\mathcal{B}^{p}_{i}(\xi)\omega_{i}$.
A NURBS curve is a weighted combination of NURBS functions defined by 
\begin{equation}\label{paramTOphys1D}
\mathcal{C}^{p}(\xi) = \sum\limits_{i=1}^{m_{b}}\mathcal{\widetilde{N}}^{i}_{p}(\xi){B}_{i},
\end{equation} 
where ${B}_{i} \in \mathbb{R}^{d}$, $i=1,\dots,m_{b}$, are called control points. In the same fashion, a NURBS surface is defined by
\begin{equation}\label{paramTOphys2D}
\mathcal{S}^{p,q}(\xi,\eta)=\sum\limits_{i=1}^{m_{b}}\sum\limits_{j=1}^{l_{b}}\mathcal{\widetilde{N}}^{i}_{p}(\xi)\mathcal{\widetilde{N}}^{j}_{q}(\eta)B_{i,j},
\end{equation}
where $B_{i,j} \in \mathbb{R}^{d}$, with $i=1,\dots,m_{b}$ and $j=1,\dots,l_{b}$, are the control points. Here, $m_{b}$ and $l_{b}$ denote the total number of control points in $x$- and $y$-directions, respectively. For simplicity, we will use the compact form of equation \eqref{paramTOphys2D} in the rest of this paper as 
\begin{equation}\label{paramTOphys2Dcompact}
\mathcal{S}^{p,q}(\xi,\eta)=\sum\limits_{m=1}^{\text{Ndof}}\mathcal{\widetilde{N}}^{m}_{p,q}(\xi,\eta)B_{m},
\end{equation} 
where $\mathcal{\widetilde{N}}^{m}_{p,q}$ is the basis function defined as $\mathcal{\widetilde{N}}^{m}_{p,q}(\xi,\eta) = \mathcal{\widetilde{N}}_{p}^{i}(\xi)\mathcal{\widetilde{N}}_{q}^{j}(\eta)$ and $\text{Ndof} = m_{b}\times l_{b}$ is the total number of control points. Note that in this compact form, the index $m$ is linked to the indices $i$ and $j$ by $m = i + (j - 1)m_{b}$.\\
\begin{figure}[t]
\begin{center}
\includegraphics[width=0.8\textwidth]{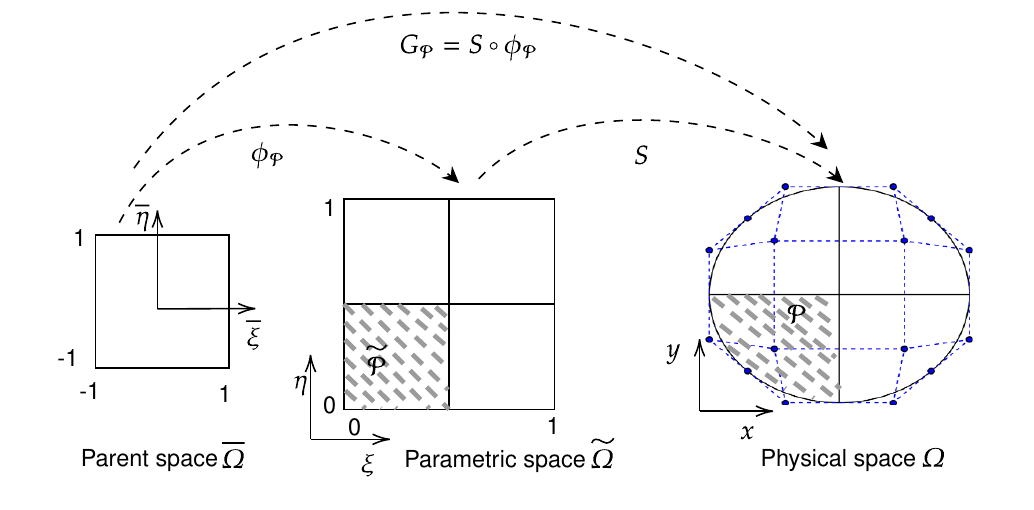}%
\vspace{-5mm}
\caption{Mappings from the parent space through the parametric space to the physical space. For the purpose of illustration, we took the knot vectors $\Xi^{1} = \Xi^{2} = \bigl\{0,0,0,1/2,1,1,1\bigr\}$. Thus, the parametric space consists of four elements. We chose the degrees $p = q = 2$ which results in $16=4\times 4$ NURBS functions obtained as tensor product of the corresponding one-dimensional NURBS functions. The dashed blue line in the physical space refers to the control polygon associated to the circular object $\Omega$.}\label{mapping}
\end{center}
\end{figure}

It should be pointed out that each element is mapped from the physical space to the parent space, also called the master space, where all IgA computations are performed, following the same procedure carried out in the conventional finite element analysis, see for instance \cite{nguyen2015isogeometric} and further details are therein. However, an additional space, which is absent in the conventional finite element analysis, is required in order to construct the NURBS basis functions which is called the parametric space. An additional mapping is required in order to operate  in the coordinates of the parent space. The mentioned mappings are illustrated in Figure \ref{mapping}, where $\phi_{\mathcal{P}}:\overline{\Omega}\longrightarrow \widetilde{\mathcal{P}}\subset \widetilde\Omega$ denotes the mapping from the parent space to the parametric space, and $S:\widetilde{\Omega} \longrightarrow \Omega$ denotes the mapping from the parametric space to the physical space. Note that this mapping satisfies $S(\widetilde{\mathcal{P}})= \mathcal{P}$. Therefore, the mapping $G_{\mathcal{P}}:\overline{\Omega}\longrightarrow \mathcal{P}$ is given by the composition $S \circ\phi_{\mathcal{P}}$. The function $\phi_{\mathcal{P}}$ maps the square $[-1,1]\times[-1,1]$ to the element $\widetilde{\mathcal{P}}=[\xi_{i},\xi_{i+1}]\times[\eta_{j},\eta_{j+1}]$, where $\xi_{i} < \xi_{i+1}$ and $\eta_{i} < \eta_{i+1}$, according to the following formula
\begin{equation}
\label{mapp}
\phi_{\mathcal{P}}(\overline{\xi},\overline{\eta}) = \left(\tfrac{1}{2}\Bigl((\xi_{i+1}-\xi_{i})\overline{\xi} + (\xi_{i+1}+\xi_{i})\Bigr),\tfrac{1}{2}\Bigl((\eta_{j+1}-\eta_{j})\overline{\eta} + (\eta_{j+1}+\eta_{j})\Bigr)\right)^T.
\end{equation}
Note that the mapping $S$ in a two-dimensional problem is given by \eqref{paramTOphys2D}. Finally, $G_{\mathcal{P}}$ can be straightforwardly written as follows
\begin{equation}
{G_{\mathcal{P}}^{p,q}(\overline{\xi},\overline{\eta})} = \sum_{m=1}^{\text{Ndof}}\mathcal{\widetilde{N}}^{m}_{p,q}\left(\phi_{\mathcal{P}}\left(\overline{\xi},\overline{\eta}\right)\right)B_{m}.
\end{equation}
It should be noted that the above approach is required to discretize space in the new IgA operator splitting proposed in this work. For the sake of brevity, we have kept the above presentation short. Further details can be found in \cite{hughes2005isogeometric,cottrell2009isogeometric,nguyen2015isogeometric,piegl2012nurbs}
among others.
\section{Isogeometric Strang splitting}\label{sec3}
In this section, we present the operator splitting technique used to approximate problem \eqref{model}. The method consists of splitting the problem into two parts, the nonlinear reaction part and the diffusion part, where the advective terms are treated using the semi-Lagrangian approach, and then solving each part using an appropriate method. It is worth mentioning that operator splitting schemes are very suitable for problems involving multiscale effects such as interaction of fluid and structure. Indeed, solving each subproblem with an appropriate technique reduces the complexity resulting from the solution of large scale problems in engineering and science. In the current work, the problem consists of an advection-diffusion-reaction system, where the reaction terms are nonlinear and sometimes involve stiffness. Moreover, the dominance of advective terms complicates the problem and leads to boundary layers with very strong gradients. Stabilization techniques such as SUPG \cite{akin2004calculation,brooks1982streamline} and SGS \cite{hauke2002simple} are often used to overcome these challenges, which may lead to additional computational cost. In the current study, we adopt a temporal divide-and-conquer procedure and solve each subproblem using an appropriate approach.\\

In order to derive the splitting used in this study, we consider first the total derivative, also called the material or substantial derivative. It measures the rate of change of the quantities $u$ and $v$ along the trajectories of flow particles. The material derivative is given as
\begin{subequations}
\begin{equation}
\label{conv-a}
\frac{D u}{D t} = \frac{\partial u}{\partial t} + {\bf a}\cdot\nabla u,
\end{equation}
\begin{equation}
\label{conv-b}
\frac{D v}{D t} = \frac{\partial v}{\partial t} + {\bf a}\cdot\nabla v.
\end{equation}
\end{subequations}
Rewriting the advective terms in terms of material derivatives, we are left with the following problem
\begin{equation}\label{StrgDif}
\begin{aligned}
\frac{D u}{D t}  - \nabla \cdot (d_{1}\nabla u) &= f(u,v), \quad \text{in} \quad \Omega \times [0,T],\\
\frac{D v}{D t}  - \nabla \cdot (d_{2}\nabla v) &= g(u,v), \quad \text{in} \quad \Omega \times [0,T].
\end{aligned}
\end{equation}
Setting $U = (u,v)^T$, problem \eqref{StrgDif} can be rewritten as
\begin{equation}\label{eq:pde-compact}
\frac{D U}{D t} = \mathcal L_1(U) + \mathcal L_2(U),
\end{equation}
where $\mathcal L_1$ denotes the nonlinear reaction operator
$$
\mathcal L_1(U) = \begin{bmatrix}  f\left(u,v\right)\\[0.7ex]
g\left(u,v\right)\end{bmatrix},
$$
and $\mathcal L_2$ denotes the diffusion operator
$$
\mathcal L_2(U) = \begin{bmatrix} - \nabla\cdot(d_{1}\nabla u) \\[0.7ex]
- \nabla\cdot(d_{2}\nabla v) \end{bmatrix}.
$$
For the time integration of \eqref{eq:pde-compact}, we use the classical Strang splitting technique which is proven to be second-order accurate, see \cite{dellar2013interpretation,lanser1999analysis,strang1968construction} for instance. For this purpose, let $\Delta t$ denote the time step size and divide the integration interval $[0,T]$ into subintervals of equal length $[t^n, t^{n+1}] = [t^n, t^{n+\frac12}]\cup [t^{n+\frac12}, t^{n+1}]$ with $t^n = n\Delta t$. Further, let $U^n(\bfx)$ denote the numerical solution at time $t^n$. Then, solving the following subproblems
\begin{subequations}
\begin{align}
\label{eqq-a}
 \frac{\partial U^*}{\partial t} &= \mathcal L_1(U^*),\quad \text{with} \quad U^{*}(t^{n},{\bfx}) = U^n({\bfx})\quad t \in [t^{n},t^{n+\frac{1}{2}}],\\
\label{eqq-b}
 \frac{D U^{**}}{D t} &= \mathcal L_2(U^{**}), \quad \text{with} \quad U^{**}(t^{n},{\bfx}) = U^{*}(t^{n+\frac{1}{2}},{\bfx}),\quad t \in [t^{n},t^{n+1}],\\
\label{eqq-c}
 \frac{\partial U^{***}}{\partial t} &= \mathcal L_1(U^{***}), \quad \text{with} \quad U^{***}(t^{n+\frac{1}{2}},{\bfx}) = U^{**}(t^{n+1},{\bfx}),\quad t \in [t^{n+\frac{1}{2}}, t^{n+1}],
\end{align}
\end{subequations}
provides the numerical solution $U^{n+1}({\bfx}) = U^{***}(t^{n+1},{\bfx})$ at time $t^{n+1}$. It should be stressed that solving the subproblem \eqref{eqq-b} requires the solution of linear systems, involving mass and stiffness matrices, which is computationally much more expensive compared to \eqref{eqq-a} and \eqref{eqq-c}. Thus, the reaction part is solved twice along the half time increments $[t^{n},t^{n+ \frac{1}{2}}]$ and $[t^{n+ \frac{1}{2}},t^{n+1}]$ in the Strang splitting, whereas the diffusion part is solved only once along the full time increment $[t^{n},t^{n+1}]$.
\subsection{Semi-Lagrangian isogeometric analysis for the advection-diffusion terms}\label{subsec3.1}
In this section we deal with the advection-diffusion terms. One possibility would be to use a purely Lagrange--Galerkin method to approximate the material time derivatives. In that approach, the particles are followed forward in time and the mesh can become distorted to such an extend that remeshing becomes necessary. In such scenarios, remeshing is expensive and can often lead to inaccuracies in the solution. In the semi-Lagrangian approach, however, the particles are tracked backward. Hence, the mesh remains fixed and does not move during time evolution which is the main advantage of semi-Lagrangian methods over their purely Lagrangian counterparts. This also reflects the strength of the semi-Lagrangian schemes which offer more convenience in accurately studying the interaction of substances $u$ and $v$ in the presence of transport phenomena. In the first step, we thus solve the advective terms by means of the semi-Lagrangian method in the framework of $L^{2}$-projection. The obtained parabolic problems are integrated in the second step. It should be stressed that one of the advantages of the NURBS-based semi-Lagrangian method, in the absence of strong-coupling such as cross-diffusion, is the decoupling of the equations for the unknowns $u$ and $v$. This allows one to compute each of them separately and immediately reduces the computational cost. Henceforth, the method is presented only for the component $u$ and the same procedure is carried out straightforwardly for the component $v$.
\subsubsection{Calculation of the departure points}
\begin{figure}[t]
\begin{center}
\includegraphics[width=0.9\linewidth]{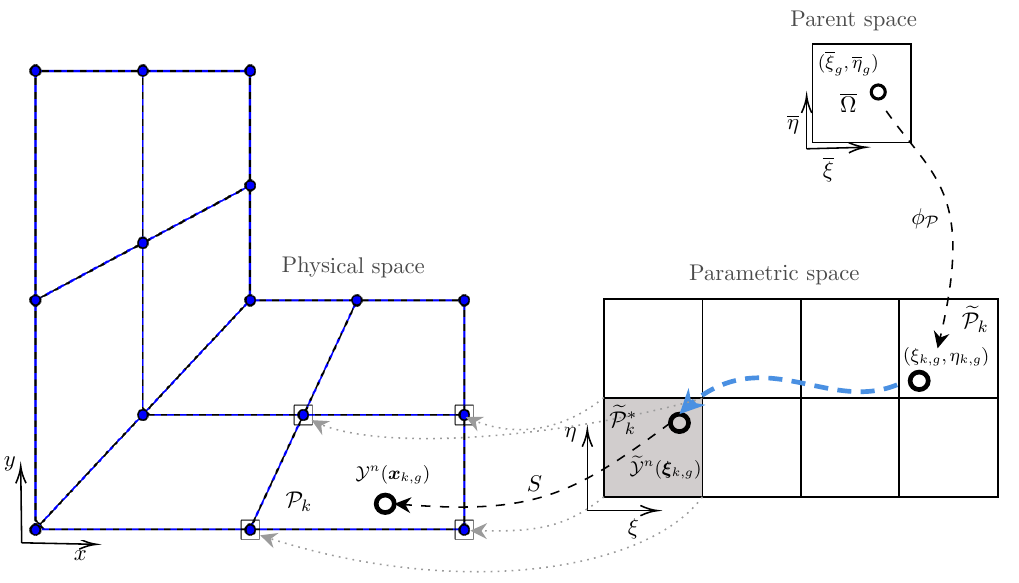}\hspace{2ex}
\caption{A schematic diagram showing the main quantities used in the calculation of the departure points. Since NURBS functions lie in the parametric space, each quadrature point $(\bar{\xi_{g}},\bar{\eta_{g}})$ is first mapped from the parent space to the parametric space according to the mapping $\phi_{\mathcal{P}_{k}}$. This results in the point $\bm\xi_{k,g}=(\xi_{k,g},\eta_{k,g})$. The corresponding departure point $\mathcal{\widetilde{Y}}^{n}(\bm\xi_{k,g})$ is then calculated at the host element, see the dark element $\widetilde{\mathcal{P}}^{*}_{k}$, before being mapped to the physical element $\mathcal{P}_{k}$. This results in the departure point $\mathcal{Y}^{n}({\bfx}_{k,g})$. In this illustration, the local control points for NURBS functions (of degree~1) are indicated by small squares.}\label{depp}
\end{center}
\end{figure}

Consider the subinterval $[t^{n},t^{n+1}]$ and let $\Delta t = t^{n+1}-t^n$ denote the time step size. The spatial discretization is done with NURBS functions. Thus, the computational domain $\Omega$ is partitioned into a set of elements $\mathcal{P}_{k}$ such that $ \Omega = \bigcup_{k=1}^{N_{e}}\mathcal{P}_{k} $ where $N_{e}$ denotes the total number of elements. We further provide the parent space with quadrature points $(\bar{\xi}_{g},\bar{\eta}_{g})$, $g=1,\dots,N_{Q}$, where $N_{Q}$ denotes the total number of such points in an element. These quadrature points are mapped to each parametric element $\widetilde{\mathcal{P}}_{k}$ using the map $\phi_{\mathcal{P}}$ and result in the points $\bm\xi_{k,g} = (\xi_{k,g},\eta_{k,g})^T$. For the $L^{2}$-projection approach which is discussed below, we will need the departure point of each quadrature point $\bm\xi_{k,g}$. These are obtained by solving the following ODEs, see Figure \ref{depp} for an illustration:
\begin{equation}\label{conv13}
\begin{aligned}
\frac{d \widetilde{\mathcal{Y}}(\bm\xi_{k,g},\tau)}{d\tau} &= {\bf a}\left(\widetilde{\mathcal{Y}}(\bm\xi_{k,g},\tau),\tau \right),\qquad  t^{n} \leq \tau \leq t^{n+1},\\
\widetilde{\mathcal{Y}}(\bm\xi_{k,g},t^{n+1}) &= \bm\xi_{k,g}.
\end{aligned}
\end{equation}
Here, $\widetilde{\mathcal{Y}}(\bm\xi_{k,g},\tau)$ denotes the departure point at time $\tau$ of a particle that will reach the quadrature point $\bm\xi_{k,g} = (\xi_{k,g},\eta_{k,g})^T$ at time $t^{n+1}$. \\

To approximate \eqref{conv13}, many schemes can be used. However, the accuracy of a semi-Lagrangian method is essentially related to the precision of the algorithm used to compute the departure points. To approximate $\widetilde{\mathcal{Y}}(\bm\xi_{k,g},\tau)$, some authors adopt a second-order explicit Runge--Kutta scheme, which may be not accurate enough to maintain a particle on its curved trajectory. In this study, we use a third-order Runge--Kutta scheme for greater accuracy and stability. Starting at the final value $\bm\xi_{k,g}$, the scheme reads
\begin{equation}\label{RK3}
\begin{aligned}
\mathcal{K}^{(1)}_{k,g} &= \bm\xi_{k,g} - \Delta t\,{\bf a}\left(\bm\xi_{k,g},t^{n+1}\right),\\
\mathcal{K}^{(2)}_{k,g} &= \frac34\bm\xi_{k,g} +\frac14 \mathcal{K}^{(1)}_{k,g} - \frac14\Delta t\,{\bf a}\left(\mathcal{K}^{(1)}_{k,g},t^{n+\frac12}\right),\\
\widetilde{\mathcal{Y}}{^n}(\bm\xi_{k,g}) &= \frac13\bm\xi_{k,g} +\frac23\mathcal{K}^{(2)}_{k,g} - \frac23\Delta t\,{\bf a}\left(\mathcal{K}^{(2)}_{k,g},t^n\right).
\end{aligned}
\end{equation}
The departure point $\widetilde{\mathcal{Y}}{^n}(\bm\xi_{k,g})$ is then mapped to the physical space through the following equation
\begin{equation}\label{RK3a}
\mathcal{Y}{^n}(\bfx_{k,g}) = \sum_{m=1}^{\text{Ndof}}\widetilde{\mathcal{N}}_{p,q}^{m}\left(\widetilde{\mathcal{Y}}{^n}(\bm\xi_{k,g})\right)B_{m}.
\end{equation}
\subsubsection{The $L^{2}$-projection BDF2 approach}
The second step of Strang splitting deals with the diffusive terms. Here, the advective terms are treated with a semi-Lagrangian approach and we are left with a parabolic problem, where the time derivatives are expressed in terms of material derivatives. As before, the method is presented only for the component $u$ and the same procedure is straightforwardly carried out for the component $v$. Following \eqref{eqq-b} and applying the second-order BDF2 with step size $\Delta t$, the weak formulation of the diffusive term related to the activator $u$ reads as follows: 
\begin{equation}\label{modelOp}
\int_{\Omega} \left( \frac{3}{2\Delta t} u_{h}^{n+1}(\bfx) - \frac{2}{\Delta t} u_{h}^{n}\left(\mathcal{Y}^{n}_{n+1}(\bfx)\right) + \frac{1}{2\Delta t} u_{h}^{n-1}\left(\mathcal{Y}^{n-1}_{n+1}(\bfx)\right) \right) \phi_{h}(\bfx)\;d\Omega  = - \int_{\Omega}d_{1}\nabla u_{h}^{n+1}(\bfx) \nabla \phi_{h}(\bfx)\;d\Omega, \quad \forall \phi_{h} \in \mathcal{V}_{h},
\end{equation}
were the set $\mathcal{V}_{h}$ denotes the conforming finite element space given by
\[
\mathcal{V}_{h} = \bigg\{{w}_{h}: \Omega \rightarrow \mathbb{R}^{2} \;\Big|\ w_{h} \circ S \in \mathcal{S}_{m}^{p,q},\; w_{h} = 0 \; \text{on} \; \partial \Omega \bigg\}.
\]
Here $\mathcal{S}_{m}^{p,q} = \mathcal{S}_{m}^{p,q}(\Xi_{1},\Xi_{2}) = \text{span}\big\{\widetilde{\mathcal{N}}_{p,q}^{m}\big\}_{m=1}^{\text{Ndof}}$ denotes the two-dimensional space of NURBS functions of degree $(p,q)$ determined by the tensor product of two knot vectors $\Xi_{1}$ and $\Xi_{2}$. Using equation \eqref{RK3a}, we calculate $\mathcal{Y}^{n}_{n+1}(\bfx)$ and $\mathcal{Y}^{n-1}_{n+1}(\bfx)$ which denote the departure points that will reach the point $\bfx$ at time $t_{n+1}$ starting from times $t_{n}$ and $t_{n-1}$, respectively. Therefore $u_{h}^{n}$ and $u_{h}^{n-1}$ are the solutions evaluated one time step and two time steps back along the characteristics.\\

The semi-discrete solution $u_{h}$ can be rewritten in terms of NURBS functions as 
\begin{equation}\label{Sbasis}
u_{h}^{n+1}({\bfx}) = \sum \limits_{l=1}^{\text{Ndof}}\bar{u}_{l}^{n+1}\mathcal{N}^{l}_{p,q}(\bfx),
\end{equation}
where $\mathcal{N}^{l}_{p,q}(\bfx)= \widetilde{\mathcal{N}}^{l}_{p,q}(S^{-1}(\bfx))$ and $\bar{u}_{l}$ is the set of the (unknown) semi-discrete control variables. 
Inserting \eqref{Sbasis} into \eqref{modelOp} and rearranging all the terms, we obtain the following compact form
\begin{equation}\label{matrixForm} 
\frac{3}{2\Delta t}{\bf M}\bar{\bf u}^{n+1} + {\bf K }\bar{\bf u}^{n+1} = \frac{2}{\Delta t} \widetilde{{\bf H}}^{n} - \frac{1}{2\Delta t}\widehat{{\bf H}}^{n-1},
\end{equation} 
where $\bar{\bf u}^{n+1}$ is the vector of the control variables $u_{m}^{n+1}$ at time $t_{n+1}$, $\bf M$ and $\bf K$ are the mass and stiffness matrices which entries $\mathcal{M}_{ml}$ and $\mathcal{K}_{ml}$ are the Gauss--Legendre approximation as follows:
\begin{equation}\label{026}
\begin{aligned}
\mathcal{M}_{ml}&= \int_{\Omega}\mathcal{N}^{l}_{p,q}(\bfx)\mathcal{N}^{m}_{p,q}(\bfx)\;d\Omega = \sum_{k=1}^{N_{e}}\int_{\mathcal{P}_{k}}\mathcal{N}^{l}_{p,q}(\bfx)\mathcal{N}^{m}_{p,q}(\bfx)\;d \Omega \\
&= \sum_{k=1}^{N_{e}}\int_{\bar{\Omega}}\widetilde{\mathcal{N}}^{l}_{p,q}\left( \phi_{{\mathcal{P}_k}}(\bar{\xi},\bar{\eta}) \right) 
\widetilde{\mathcal{N}}^{m}_{p,q}\left( \phi_{{\mathcal{P}_k}}(\bar{\xi},\bar{\eta})\right)\left|J_k(\bar{\xi},\bar{\eta})\right|\;d\bar{\Omega}\\
&\approx \sum_{k=1}^{N_{e}}\sum_{g=1}^{N_{Q}}w_{k,g}\widetilde{\mathcal{N}}^{l}_{p,q}\left( \phi_{{\mathcal{P}_k}}(\bar{\xi}_{g},\bar{\eta}_{g})\right)
\widetilde{\mathcal{N}}^{m}_{p,q}\left( \phi_{{\mathcal{P}_k}}(\bar{\xi}_{g},\bar{\eta}_{g})\right)\left|J_k(\bar{\xi}_{g},\bar{\eta}_{g})\right|
\end{aligned}
\end{equation}
and
\begin{equation}\label{027}
\begin{aligned}
\mathcal{K}_{ml}&= \int_{\Omega}d_{1}\nabla \mathcal{N}^{l}_{p,q}(\bfx)\cdot\nabla \mathcal{N}^{m}_{p,q}(\bfx)\;d\Omega = \sum_{k=1}^{N_{e}}\int_{\mathcal{P}_{k}}d_{1}\nabla \mathcal{N}^{l}_{p,q}(\bfx)\cdot\nabla \mathcal{N}^{m}_{p,q}(\bfx)\;d\Omega, \\
&= \sum_{k=1}^{N_{e}}\int_{\bar{\Omega}}d_{1}\nabla \widetilde{\mathcal{N}}^{l}_{p,q}\left( \phi_{{\mathcal{P}_k}}(\bar{\xi},\bar{\eta})\right)\cdot\nabla 
\widetilde{\mathcal{N}}^{m}_{p,q}\left( \phi_{{\mathcal{P}_k}}(\bar{\xi},\bar{\eta})\right)\left|J_k(\bar{\xi},\bar{\eta})\right|\;d\bar{\Omega},\\
&\approx \sum_{k=1}^{N_{e}}\sum_{g=1}^{N_{Q}}w_{k,g}d_{1}\nabla \widetilde{\mathcal{N}}^{l}_{p,q}\left( \phi_{{\mathcal{P}_k}}(\bar{\xi}_{g},\bar{\eta}_{g})\right)
\cdot\nabla \widetilde{\mathcal{N}}^{m}_{p,q}\left( \phi_{{\mathcal{P}_k}}(\bar{\xi}_{g},\bar{\eta}_{g})\right)\left|J_k(\bar{\xi}_{g},\bar{\eta}_{g})\right|,
\end{aligned}
\end{equation}
where $|J_k(\bar{\xi},\bar{\eta})|$ is the determinant of the Jacobian of the map $G$ from the parent space to the element $\mathcal{P}_{k}$ in physical space.
In the same manner, the entries of the right-hand sides $\widetilde{{\bf H}}^{n}$ and $\widehat{{\bf H}}^{n-1}$ are the approximation of the following formulas
\begin{equation}\label{integRHS}
\begin{array}{cc} 
H^{n}_{m}: = \int_{\Omega}u^{n}(\mathcal{Y}_{n+1}^n(\bfx))\mathcal{N}^{m}_{p,q}(\bfx)\;d\Omega, & \quad H^{n-1}_{m}:= \int_{\Omega}u^{n-1}(\mathcal{Y}_{n+1}^{n-1}(\bfx))\mathcal{N}^{m}_{p,q}(\bfx)\;d\Omega.
\end{array}
\end{equation}
To compute the integrals in the above formula we use an $L^{2}$-projection approach since the control points do not lie on the physical geometry which makes NURBS functions non-interpolatory. Note that the exact evaluation of integrals in finite element methods is not possible, in general. Therefore, the use of quadrature rules is inevitable. These numerical integrations play a crucial rule in the efficiency and accuracy of Galerkin-based finite element methods. In addition to the isogeometric paradigm of assembling systems of equations which requires the calculation of integrals for the mass and stiffness matrices, the $L^{2}$-projection approach necessitates the calculation of integrals in \eqref{integRHS}. With the help of a Gauss--Legendre quadrature rule, the integral $H^{n}_{m}$ in \eqref{integRHS} is approximated as follows
\begin{equation}\label{24}
\begin{aligned}
\int_{\Omega}u^{n}(\mathcal{Y}^n_{n+1}(\bfx))\mathcal{N}^{m}_{p,q}(\bfx)\;d\Omega&= \sum_{k=1}^{N_{e}}\int_{\mathcal{P}_{k}}u^{n}\left( \mathcal{Y}^n_{n+1} (\bfx)\right) \mathcal{N}^{m}_{p,q}(\bfx)\;d\Omega \\
&= \sum_{k=1}^{N_{e}}\int_{\bar{\Omega}}u^{n}\left( \mathcal{Y}^n_{n+1} \left(G^{p,q}_{\mathcal{P}}(\bar{\xi},\bar{\eta})\right)\right)\widetilde{\mathcal{N}}^{m}_{p,q}\left(\phi_{{\mathcal{P}}_k}(\bar{\xi},\bar{\eta})\right)\left|J_k(\bar{\xi},\bar{\eta})\right|\;d\bar{\Omega}, \\
&\approx \sum_{k=1}^{N_{e}}\sum_{g=1}^{N_{Q}}w_{k,g}u^{n}\left( \mathcal{Y}^n_{n+1} \left(G^{p,q}_{\mathcal{P}}(\bar{\xi}_{g},\bar{\eta}_{g})\right)\right)\widetilde{\mathcal{N}}^{m}_{p,q}\left(\phi_{{\mathcal{P}}_k}(\bar{\xi}_{g},\bar{\eta}_{g})\right)\left|J_k(\bar{\xi}_{g},\bar{\eta}_{g})\right|, \\
&= \sum_{k=1}^{N_{e}}\sum_{g=1}^{N_{Q}}w_{k,g}\widetilde{u}_{k,g}^{n}\widetilde{\mathcal{N}}^{m}_{p,q}\left(\phi_{{\mathcal{P}}_k}(\bar{\xi}_{g},\bar{\eta}_{g})\right)\left|J_k(\bar{\xi}_{g},\bar{\eta}_{g})\right| = \widetilde{H}_m^{n}.
\end{aligned}
\end{equation}
Here, $\widetilde{u}_{k,g}^{n} = u^{n}(\mathcal{Y}^n_{n+1}(\bfx_{k,g}))$ is the solution at the departure point $\mathcal{Y}^n_{n+1}(\bfx_{k,g})$ for which the values are calculated as 
\begin{equation}\label{isoSOL}
\widetilde{u}_{k,g}^{n} = u^{n}\left(\mathcal{Y}^n_{n+1}(\bfx_{k,g})\right) \approx \sum_{l=1}^{\text{Ndof}}\bar{u}^{n}_{l}\mathcal{N}^{l}_{p,q}\left(\mathcal{Y}^n_{n+1} (\bfx_{k,g})\right),
\end{equation}
where $\bar{u}^{n}_{l}$ denote the known solutions at the control points. Further note that $w_{k,g}$ denote the quadrature weights of the Gauss--Legendre quadrature rule used in our study to evaluate the integrals. The integral $H^{n-1}_{m}$ is approximated in the same way as equation \eqref{24} where the solution in \eqref{isoSOL} is evaluated two time steps back along the characteristics using $\mathcal{Y}^{n-1}_{n+1}(\bfx_{k,g})$.\\

For the practical implementation of \eqref{isoSOL}, we take advantage of the fact that the NURBS functions have compact support. To determine the potentially non-zero elements of \eqref{isoSOL}, we perform the following search algorithm. Given $k$, we determine the element $\mathcal{P}^{*}_{k}$ to which the departure point $\mathcal{Y}^n_{n+1}(\bfx_{k,g})$ belongs. Next, we identify all indices $l$ such that $\text{supp}\,\Big( \widetilde{\mathcal{N}}^{m}_{p,q}\Big) \cap \mathcal{P}^{*}_{k} \not= \emptyset$ and use only those indices in \eqref{isoSOL}, see Figure \ref{depp} for an illustration. This local strategy greatly reduces the computational cost.\\

System \eqref{matrixForm} can be rewritten in the following more compact form
\begin{equation}\label{compaForm} 
{\bf A} \bar{\bf u}^{n+1} = {\bf B},
\end{equation} 
where ${\bf A} = \left(\frac{3}{2\Delta t} {\bf M} + {\bf K}\right)$ and ${\bf B} = \left(\frac{2}{\Delta t} \widetilde{{\bf H}}^{n} - \frac{1}{2\Delta t}\widehat{{\bf H}}^{n-1}\right)$. To solve \eqref{compaForm}, iterative solvers \cite{garcia2018refined}, overlapping domain decomposition techniques \cite{da2012overlapping,cho2020overlapping}, multilevel and multigrid techniques \cite{gahalaut2013multigrid,hofreither2017robust}, among others, can be used. Note that the condition number increases exponentially with increasing NURBS degree, so the global matrix becomes ill-conditioned and thus the iterative solvers may fail to converge in this case \cite{garcia2017value}, while direct solvers are widely used in IgA
and are known to give better results \cite{garcia2017value,garcia2017optimally,garcia2019refined}. In addition, the linear system must be evaluated and solved at each time step. Therefore, the solution process can be computationally demanding, if the system contains a large number of equations. In the present study, we use the multifrontal direct solver first introduced in \cite{duff1983multifrontal,amestoy2001fully}. There, the system matrix may be decomposed at the first time step and retained for reuse, with just updating the right-hand side of the linear system of equations at all the subsequent time steps. Standard finite element approximations often use low-order polynomial basis functions that require a large number of degrees of freedom to achieve high accuracy. This makes the solution prohibitively expensive in terms of computation. As an alternative, we propose the use of high-order basis functions, which can significantly reduce the number of unknowns required, and hence, much smaller linear systems have to be solved during the time integration procedure.

\subsection{Treatment of the reaction terms}\label{subsec3.2}
The first and the third step of Strang splitting require the solution of the following weak form of reaction equations:
\begin{eqnarray}
\int_{\Omega}\frac{\partial u_{h}}{\partial t}\phi_{h}\;d\Omega  &=& \int_{\Omega}f(u_{h},v_{h})\,\phi_{h}\;d\Omega, \quad \forall \phi_{h} \in \mathcal{V}_{h},\nonumber\\[-1ex]
\label{weakReac}\\[-1ex]
\int_{\Omega}\frac{\partial v_{h}}{\partial t}\psi_{h}\;d\Omega  &=& \int_{\Omega}g(u_{h},v_{h})\,\psi_{h}\;d\Omega, \quad \forall \psi_{h} \in \mathcal{V}_{h}.\nonumber
\end{eqnarray}
Let $U_{h} = (u_{h},v_{h})^T$, $\Phi_{h} = (\phi_{h},\psi_{h})^T$ and $F = (f,g)^T$. Then, equations \eqref{weakReac} can be rewritten in the following compact form
\begin{equation}\label{weakRsys}
\int_{\Omega}\frac{\partial U_{h}}{\partial t}\odot \Phi_{h}\;d\Omega  = \int_{\Omega}F\left(U_{h}\right)\odot \Phi_{h}\;d\Omega, \quad \forall \Phi_{h} \in \mathcal{V}_{h},
\end{equation}
where the notation $\odot$ refers to the Hadamard (or entry-wise) product, i.e., $\frac{\partial U_{h}}{\partial t}\odot \Phi_{h} = \left(\frac{\partial u_{h}}{\partial t} \phi_{h}, \frac{\partial v_{h}}{\partial t} \psi_{h}\right)^T$. The semi-discrete solution vector $U_{h}$ is rewritten in terms of NURBS functions as 
\[U_{h}(\bfx,t) = \sum \limits_{l=1}^{\text{Ndof}}\bar{U}_{l}(t)\mathcal{N}^{l}_{p,q}(\bfx).\]
The substitution of the above expressions into the system \eqref{weakRsys} yields
\begin{equation}\label{weakRsd}
\left(\sum\limits_{l=1}^{\text{Ndof}}\int_{\Omega}\mathcal{N}^{m}_{p,q}({\bfx})\mathcal{N}^{l}_{p,q}({\bfx})\; d\Omega\right)\frac{d\bar{U}_{l}^{n}}{d t}  = \int_{\Omega}F\left(\sum \limits_{l=1}^{\text{Ndof}}\bar{U}_{l}^{n}\mathcal{N}^{l}_{p,q}({\bfx})\right)\mathcal{N}^{m}_{p,q}({\bfx})\;d\Omega, \quad \forall \mathcal{N}^{m}_{p,q} \in \mathcal{V}_{h}, \quad m = 1,\dots, \text{Ndof}.
\end{equation}
The semi-discrete equations \eqref{weakRsd} can be rewritten in a compact form as 
\begin{equation}\label{weakRc}
 {\bf M}\frac{d \bar{\bf U}}{d t} = \mathcal{F}(\bar{\bf U}),
\end{equation}
where $\bf M$ denotes the mass matrix and the term $\mathcal{F}(\bar{\bf U})$ is given by the integral on the right-hand side of equation \eqref{weakRsd}. It should be stressed that this integral is approximated by the virtue of a Gauss--Legendre quadrature rule. The coupled system \eqref{weakRc} must be solved by a numerical method that is sufficiently accurate. In this work, we simply perform a step of size $\Delta t/2$ with RK4, an explicit Runge--Kutta method of order 4, see \cite{Hairer1993-vw}. For nonstiff reaction problems, this simple choice provides a sufficiently accurate solution, in general. However, for rapidly changing solutions or mildly stiff problems, one has to modify the integrator. An alternative would be to use an embedded explicit Runge--Kutta method, which automatically selects an appropriate step size to stably integrate the reaction equations. Instead of an embedded scheme, one could also use RK4 with a sub-stepping technique, e.g., performing $N_{stp}$ steps with the step size $\frac{\Delta t}{2N_{stp}}$. In order not to choose $N_{stp}$ too large, some knowledge of the stiffness is required. A standard reference for the numerical solution of ordinary differential equations are the textbooks \cite{Hairer1993-vw,wanner1996solving}.\\

Note that an implicit-explicit operator splitting is used in this study, which on the one hand does not impose serious constraints on the time steps due to the semi-Lagrangian treatment of the advective terms and the implicit integrator for the linear diffusion terms. On the other hand, the algorithm is highly accurate and easy to implement for the nonlinear reaction part since we use an explicit time discretization. By doing so, no linearization is required. Thus, given the initial solution $\bar{\bfU}^{n}$, the main steps used in the proposed operator splitting IgA method to solve the advection-diffusion-reaction problem \eqref{model} along a time step $[t^{n},t^{n+1}]$ are summarized in Algorithm \ref{Algo}.
\begin{algorithm}[H]
\caption{One time step of the isogeometric analysis Strang splitting approach for solving the advection-diffusion-reaction system \eqref{model}.}\label{Algo}
\begin{algorithmic}[1]
\STATE{Use the Runge--Kutta scheme with initial data $\bar{\bfU}^{n}$, to solve the reaction problem \eqref{weakRc} for $t^{n} \leq t \leq t^{n+\frac{1}{2}}$ to obtain $\bar{\bfU}^{n+\frac{1}{2}}$.}
\FOR{each element ${\mathcal{P}}_k$ $(k=1,\dots,N_{e})$}
\FOR{each quadrature point ${\bfx}_{k,g}$, $g=1,\dots,N_{g}$}
\STATE{Calculate the departure point $\mathcal{Y}^n_{n+1}({\bfx}_{k,g})$ according to \eqref{RK3} and \eqref{RK3a}. In the same way, compute the departure point $\mathcal{Y}^{n-1}_{n+1} ({\bfx}_{k,g})$.}
\STATE{Calculate the solution value $\widetilde{u}_{k,g}$ at $\mathcal{Y}^n_{n+1}({\bfx}_{k,g})$ using the approximation \eqref{isoSOL} with data $\bar{u}^{n+\frac{1}{2}}$ instead of $\bar{u}^{n}$. In the same way, compute $\widetilde{v}_{k,g}$. Compute also the corresponding solutions at $\mathcal{Y}^{n-1}_{n+1}({\bfx}_{k,g})$.}
\ENDFOR
\STATE{Compute the elements $\widetilde{H}_{m}^{n}$ using \eqref{24} and assemble the right-hand side vector $\widetilde{\bf H}^{n}$. Compute also the right-hand side vector $\widehat{\bf H}^{n-1}$. In the same way, compute the corresponding right-hand sides for $\widetilde{v}^{n}$.}
\ENDFOR
\STATE{Use BDF2 to solve the advection-diffusion problem \eqref{matrixForm} and the corresponding problem for $v$, for $t^{n} \leq t \leq t^{n+1}$ to obtain $\hat{\bfU}^{n+1}$.}
\STATE{Use the Runge--Kutta scheme with initial data $\hat{\bfU}^{n+1}$ to solve the reaction problem \eqref{weakRc} for $t^{n+\frac{1}{2}} \leq t \leq t^{n+1}$ to obtain $\bar{\bfU}^{n+1} = \left(\bar{\bfu}^{n+1},\bar{\bfv}^{n+1}\right)^T$.}
\end{algorithmic}
\end{algorithm}
\section{Numerical results}\label{sec4}
In this section, several numerical tests are presented to examine the performance of the new method. For the sake of brevity, all computational results are reported for only one component of the system \eqref{model}.
\begin{figure}
\begin{center}
\includegraphics[width=0.4\linewidth]{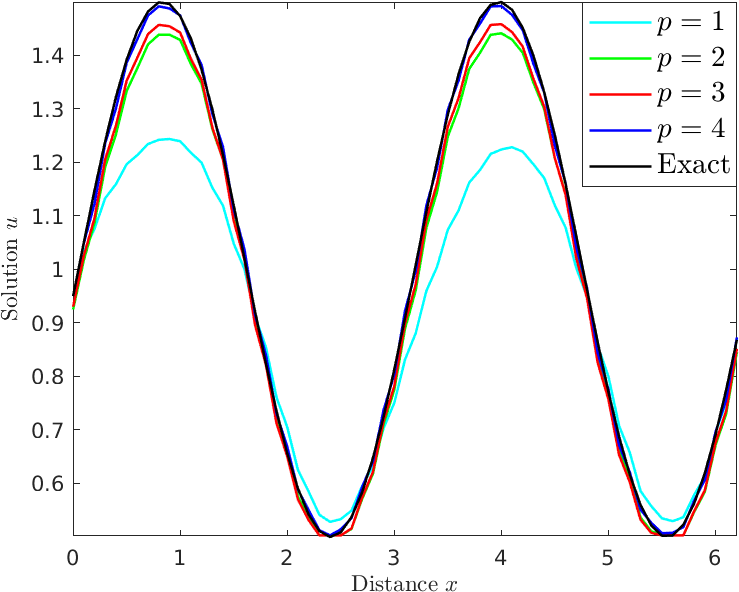}
\includegraphics[width=0.4\linewidth]{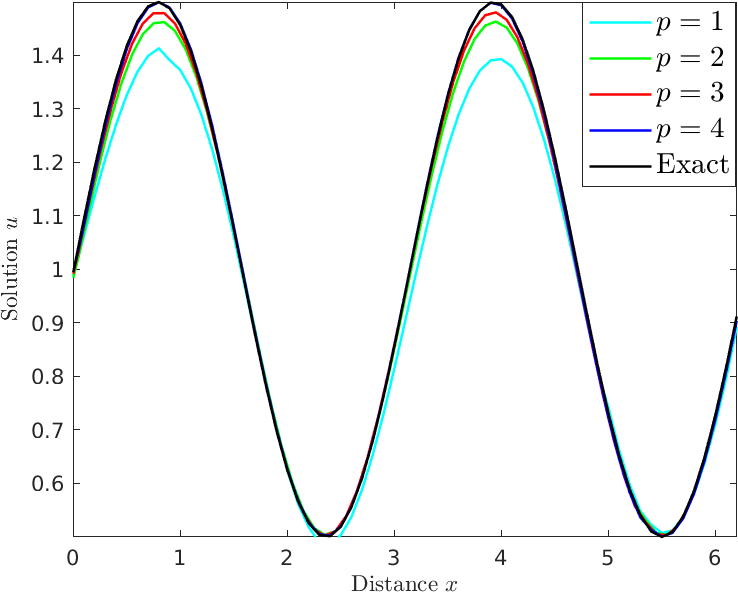}
\caption{Cross sections of the solution in Example \ref{test0} along the diagonal $x = y$, at time $t=1$ using different NURBS degrees, for two different mesh resolutions. Left: $32\times32$ elements. Right: $64\times64$ elements.}\label{Ex0_Fig1}
\end{center}
\end{figure}

\subsection{Fully nonlinear problem}\label{test0}
The first test case is examined to ascertain the accuracy of the new method for a scalar fully nonlinear equation. The problem is given by
\begin{eqnarray}\label{ex0}
\frac{\partial u}{\partial t} + u\cdot \nabla u &=&  \nabla \cdot (u \nabla u) - u^{2} + \phi({\bfx},t),\qquad \text{in} \quad \Omega \times [0,t],\nonumber\\
u({\bfx},0)  &=& u_{0}({\bfx}), \qquad \text{in} \quad \Omega,\nonumber
\end{eqnarray}
for $\Omega = [0,2\pi]^{2}$ and where the source term $\phi$ and the initial condition $u_{0}$ are calculated such that the analytical solution is given by
\[u(x,y,t) = 1 + \frac{1}{2}\sin(x + y - t).\]
Neumann boundary conditions are considered. This problem has served as a prototype to investigate the performance of several algorithms for nonlinear problems, see for example \cite{jiang2013krylov}. A notable feature of this benchmark is the need to deal with a fully nonlinear problem with the nonlinearity appearing in the advection term, the diffusion term as well as in the reaction part. It is worth noting that in this example, the stiffness matrix is assembled within the time loop since the diffusion coefficient depends on the unknown $u$. For this purpose, a Gauss--Legendre quadrature rule is used to carefully find out the support of the integrand during the assembly. The above analytical solution is used to quantify the accuracy of the proposed IgA operator splitting method by exploring different NURBS degrees. A NURBS mesh is used to discretize the computational domain $\Omega = [0,2\pi]^{2}$. This divides the domain into quadrilaterals constructed by NURBS functions. Two mesh resolutions are considered, one with element side length $h = 1/32$ and the other with element side length $h = 1/64$. To improve the accuracy of the results, refinements are performed using the $k$-refinement technique in order to ensure that the full continuity across elements is maintained \cite{hughes2005isogeometric}. In Figure \ref{Ex0_Fig1}, we show the obtained numerical solutions along the main diagonal $x = y$, using different NURBS degrees ranging from $p = 1$ to $p = 4$, on a mesh of $32\times32$ elements. The numerical solutions with the selected NURBS degrees are then repeated but with a refined mesh, as shown on the right-hand side of Figure \ref{Ex0_Fig1}. The corresponding cross sections of the exact solution are also shown for comparison reasons. The solution $u$ consists of a wave that travels continuously from one point to another in the domain and carries energy. The clear indication from Figure \ref{Ex0_Fig1} is that the linear solution (i.e., $p=1$) fails to maintain the shape of the solution $u$. Observe the large deformation that occurs in the wave along the time evolution, especially on the coarse mesh with $32\times32$ elements. The large diffusion effects that occur in the linear solution are obviously better resolved in the refined mesh of $64\times64$ elements. The quadratic and cubic solutions give roughly similar results, with small differences in the maximum value of the numerical solutions. In fact, these solutions still show some inaccuracies especially in the regions of sharp gradients. However, when the NURBS degree is increased to $p = 4$, the results are in good agreement with the analytical solution. Clearly, the shape of the wave is well preserved, without any distortion observed during the time evolution. We can clearly see that the quartic NURBS here is sufficient enough to reproduce satisfactory results for this test case with the used mesh resolution.

\begin{figure}[t]
\begin{center}
\includegraphics[width=0.23\linewidth]{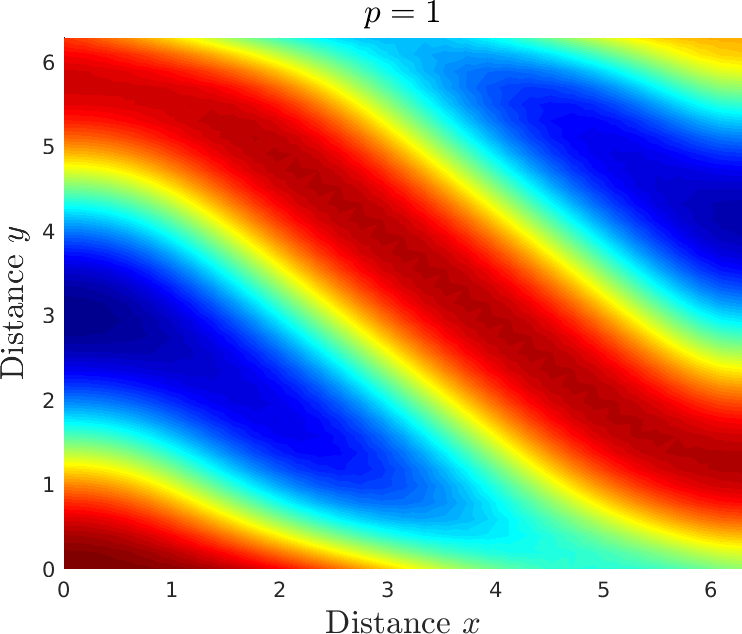}
\includegraphics[width=0.23\linewidth]{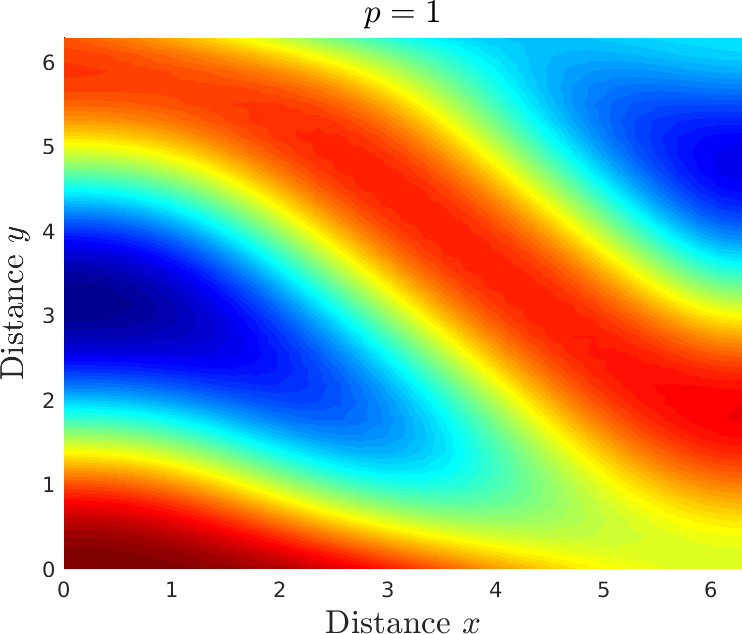}
\includegraphics[width=0.23\linewidth]{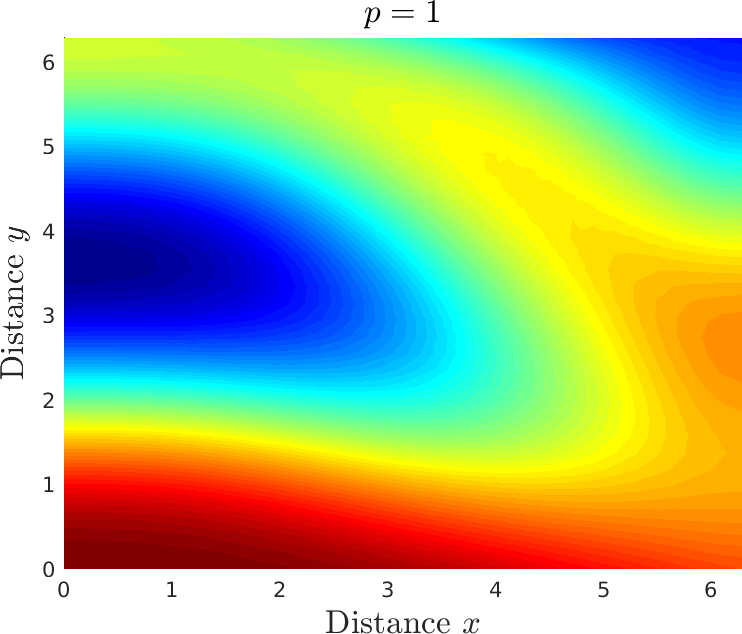}\\\vspace{0.1cm}
\includegraphics[width=0.23\linewidth]{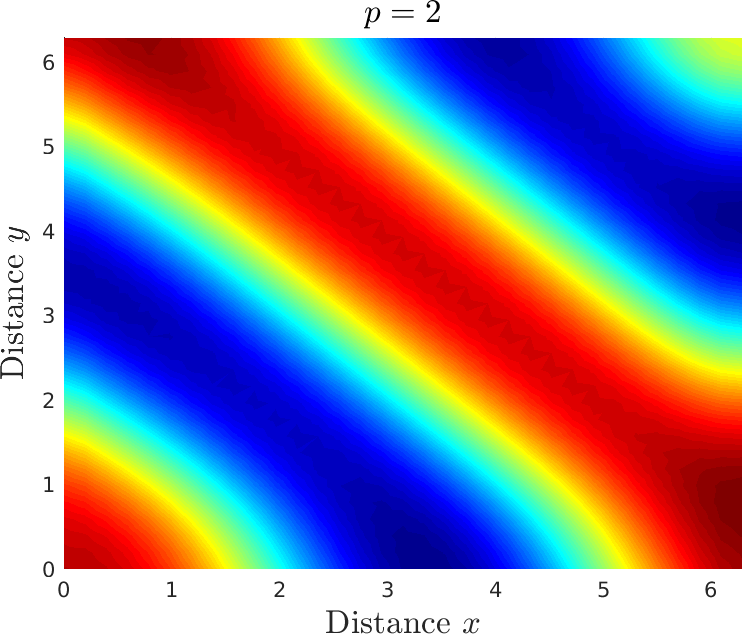}
\includegraphics[width=0.23\linewidth]{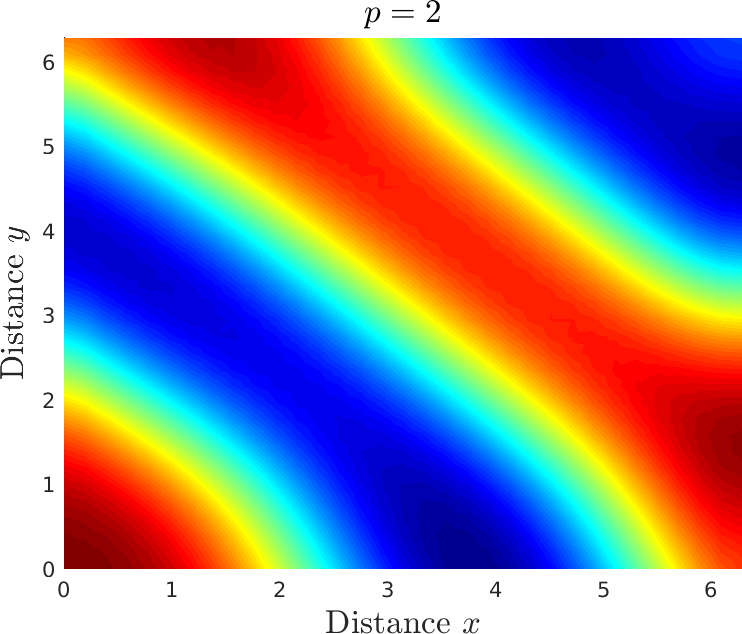}
\includegraphics[width=0.23\linewidth]{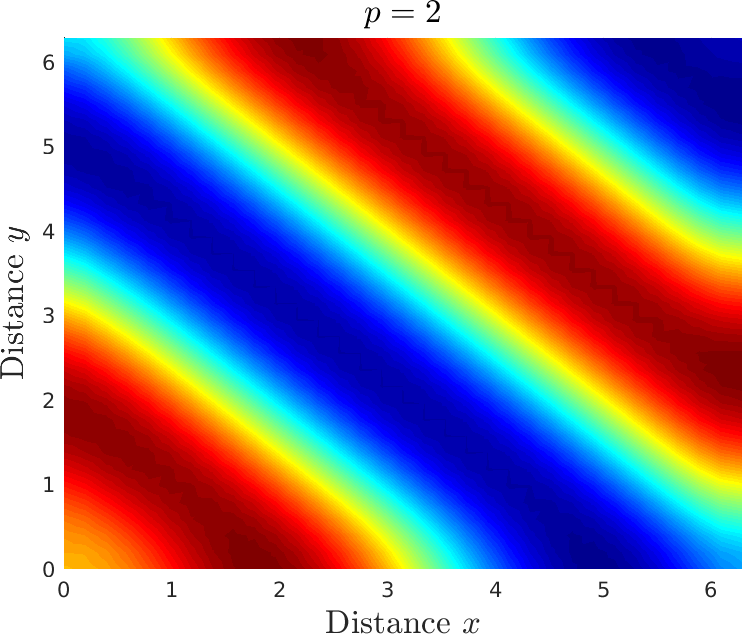}\\\vspace{0.1cm}
\includegraphics[width=0.23\linewidth]{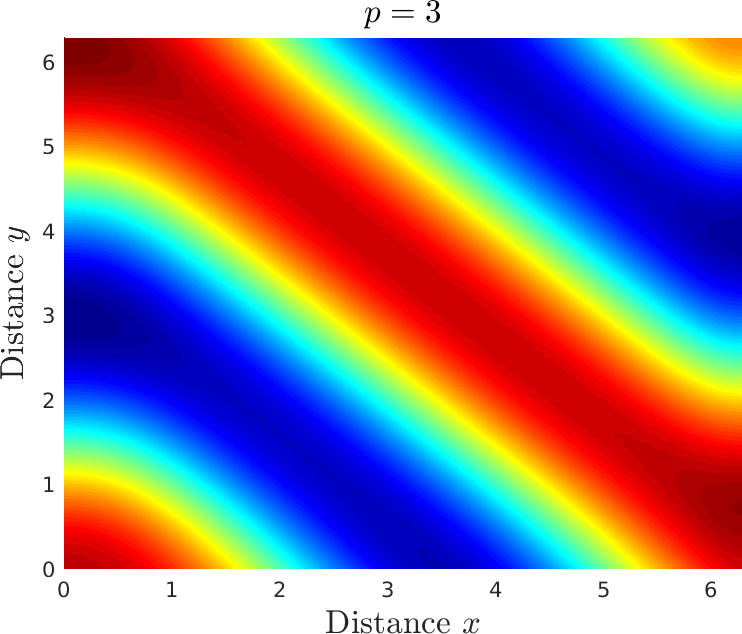}
\includegraphics[width=0.23\linewidth]{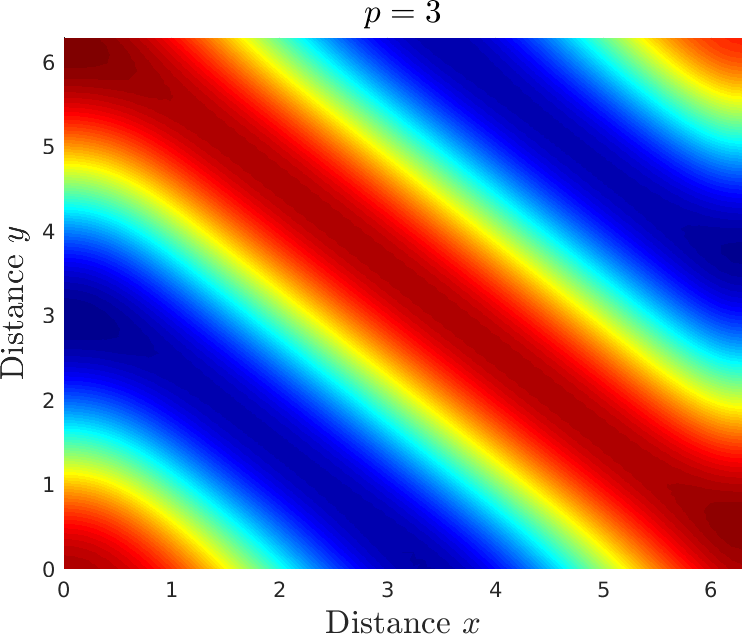}
\includegraphics[width=0.23\linewidth]{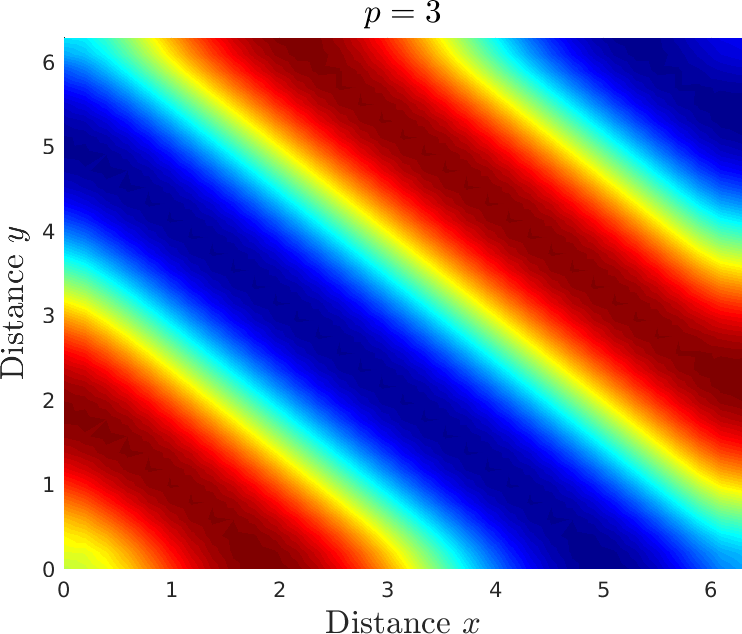}\\
\caption{Snap shots of the solution of Example \ref{test1} using different NURBS degrees. From left to right: $t = 0.5$, $t = 1$ and $t = 2$. The mesh resolution is fixed to $32\times32$ elements.}\label{Ex1_Fig1}
\end{center}
\end{figure}

\begin{figure}[h!]
\begin{center}
\includegraphics[width=0.23\linewidth]{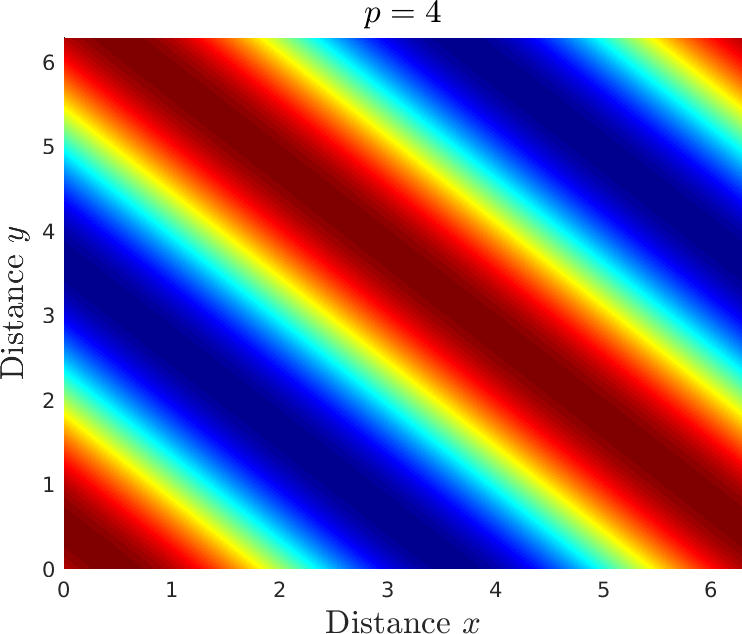}\hspace{0.1cm}\includegraphics[width=0.23\linewidth]{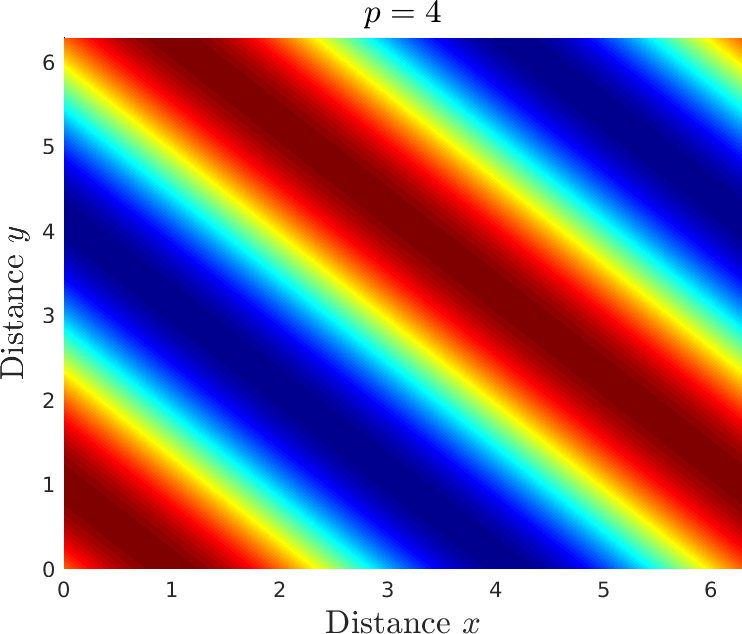}\hspace{0.1cm}
\includegraphics[width=0.23\linewidth]{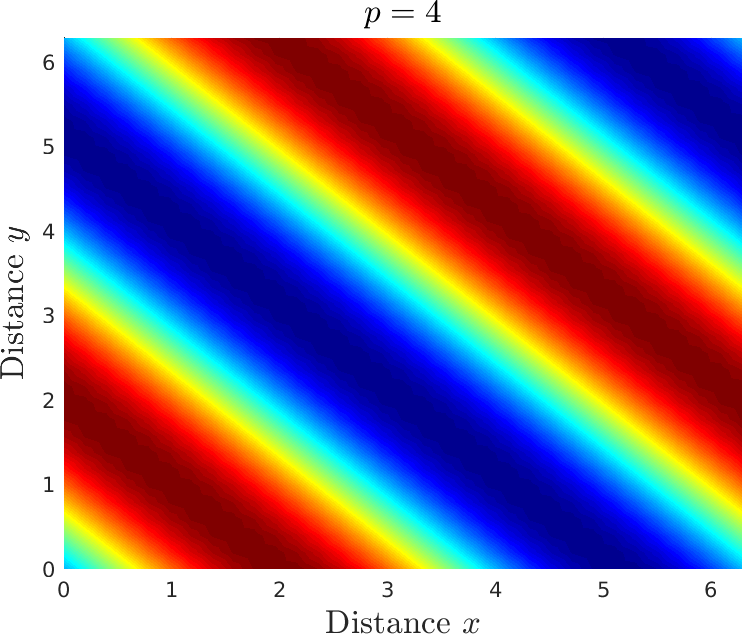}\\\vspace{0.1cm}
\includegraphics[width=0.23\linewidth]{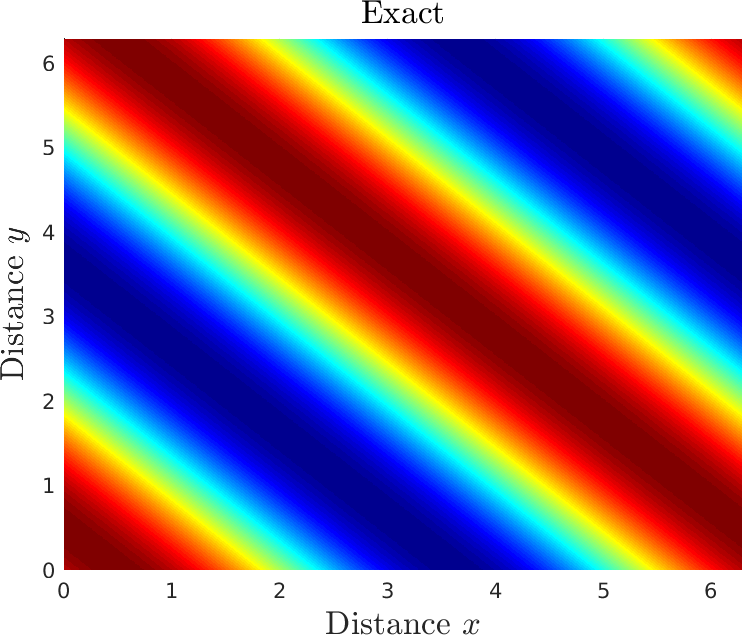}\hspace{0.06cm}
\includegraphics[width=0.23\linewidth]{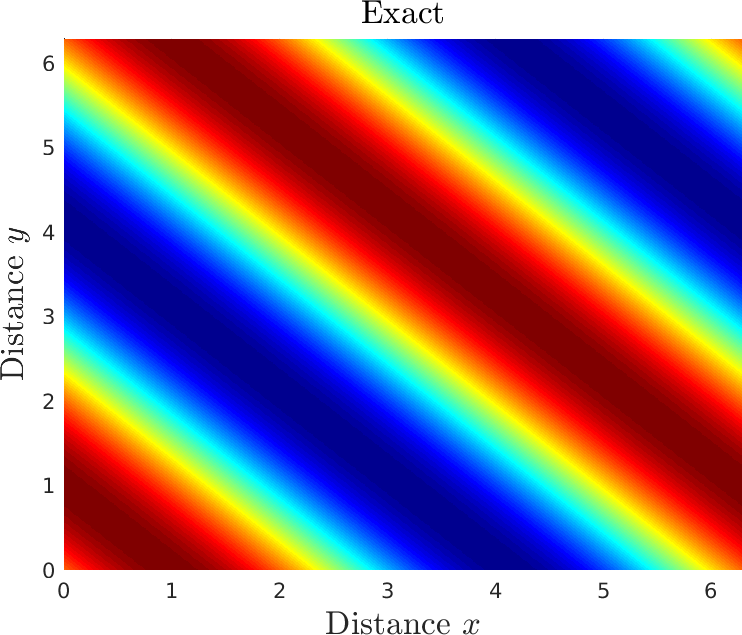}\hspace{0.06cm}
\includegraphics[width=0.23\linewidth]{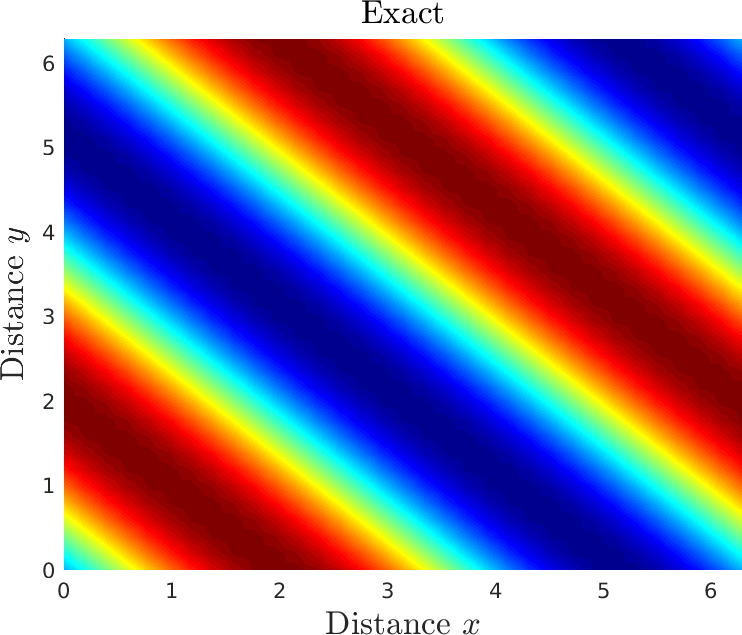}
\caption{Continuation of Fig. \ref{Ex1_Fig1}.}\label{Ex1_Fig1a}
\end{center}
\end{figure}   
\subsection{A system with exact solution}\label{test1}
We now turn our attention to studying the performance of the new method on a class of coupled equations. The system is of the form \eqref{model}, where the reaction functions are given by
\[f(u,v) = -bu - v, \quad g(u,v) = -cv,\]
with the parameters $b = 100$ and $c = 1$. The diffusion coefficients are set to $d_{1} = d_{2} = 1/2$ and the velocity field is chosen to be ${\bf a}({\bfx},t) = (a_{1},a_{2}) = (1/2,1/2)^{T}$. Boundary conditions and initial conditions are taken such that the analytical solution is given by
\begin{equation}\label{test1exact}
\begin{aligned}
u(x,y,t) &=  \big(e^{-(b + d)t} + e^{-(c + d)t}\big)\cos(x + y - at),  \\
v(x,y,t) &=  (b - c)\big(e^{-(c + d)t}\big)\cos(x + y - at),
\end{aligned}
\end{equation}
where $d = a = 1$. This problem has been studied previously in \cite{jiang2013krylov} using Krylov implicit integration factor WENO methods. Here, the analytical solution \eqref{test1exact} is used to quantify the accuracy of the proposed IgA operator splitting method using different NURBS degrees. A NURBS mesh is used to discretize the computational domain, which is given by the square $\Omega = [0,2\pi]^{2}$. The domain is divided into $32\times32$ quadrilaterals constructed with NURBS functions with element side length $h = 1/32$. To improve the accuracy of the results, refinements are performed again using the $k$-refinement technique. The obtained snapshots of the numerical solution are shown in Figures \ref{Ex1_Fig1} and \ref{Ex1_Fig1a} using different NURBS degrees ranging from $p = 1$ to $p = 4$ at three different times $t = 0.5$, $t = 1$ and $t = 2$. The corresponding snapshots of the exact solution are also shown for comparison reasons. The initial shape of the $u$ component consists of waves that propagate through the domain over time. The clear indication from Figure \ref{Ex1_Fig1} is that the linear and quadratic solutions fail to preserve the shape of the $u$ component. Note the strong deformation of the waves over time. The cubic solution still shows some deformation, especially in the boundary regions. However, when the NURBS degree is increased to $p = 4$, the results are more accurate and the solution is close to the analytical solution, see Figure \ref{Ex1_Fig1a}. In fact, the shape of the wave is well reproduced during the time evolution, and the quartic NURBS is then largely sufficient to obtain satisfactory results for this test case. \\  

\begin{figure}[b]
\begin{center}
\includegraphics[width=0.3\linewidth]{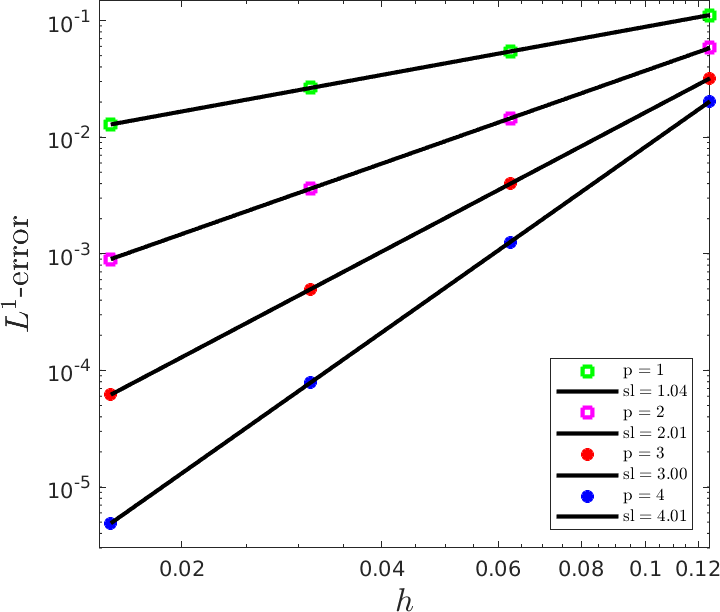}\hspace{2ex}
\includegraphics[width=0.3\linewidth]{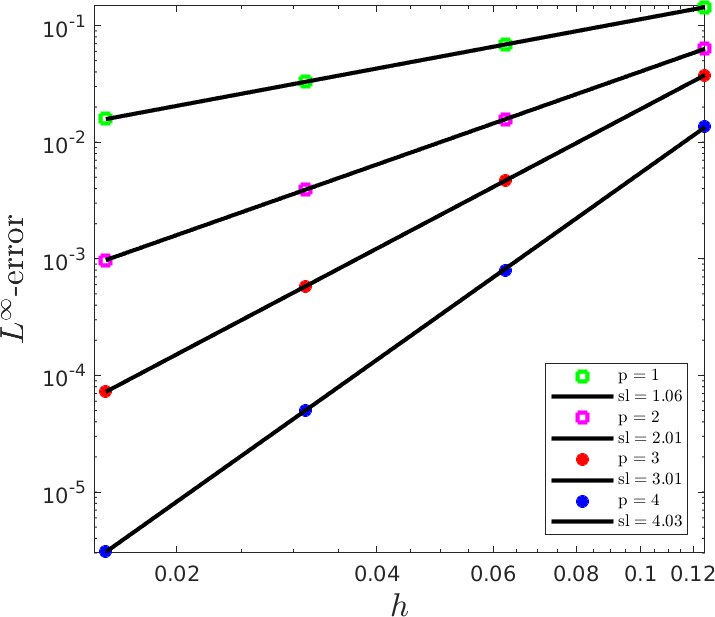}\hspace{2ex}
\caption{Convergence plots in the $L^1$ and the $L^\infty$ norm for Example \ref{test1} at time $t=1$, using different NURBS degrees. sl refers to the slop of convergence.}\label{Ex1_Fig2}
\end{center}
\end{figure}

To further determine the accuracy of the proposed IgA operator splitting method, we report in Figure \ref{Ex1_Fig2} the convergence plots in $L^1$ and $L^\infty$ norms at time $t = 1$ for different NURBS degrees. The $L^{1}$ and $L^{\infty}$ errors are plotted against the mesh size $h$. Clearly, all errors decrease as we refine the mesh and the optimum order of accuracy is achieved for all NURBS degrees. For instance, at $t = 1$ the slopes of the convergence plots for linear, quadratic, cubic and quartic NURBS are $1.06$, $2.02$, $3.01$ and $4.00$ respectively in the $L^\infty$ norm. The convergence plots in Figure \ref{Ex1_Fig2} also serve to illustrate the high accuracy achieved by the proposed IgA operator splitting method when using the $k$-refinement technique, which is an outstanding aspect of isogeometric analysis. In fact, the $k$-refinement essentially consists in increasing the smoothness across elements and eventually improves the accuracy and the efficiency compared to the classical $p$-refinement. Note, for example, that a pointwise accuracy of $10^{-4}$ is obtained for quartic NURBS with $h = 0.04$, while the same accuracy is achieved for cubic NURBS with about $h = 0.02$, which requires four times more elements for cubic NURBS. 

\begin{figure}[t]
\begin{center}
\includegraphics[width=0.23\linewidth]{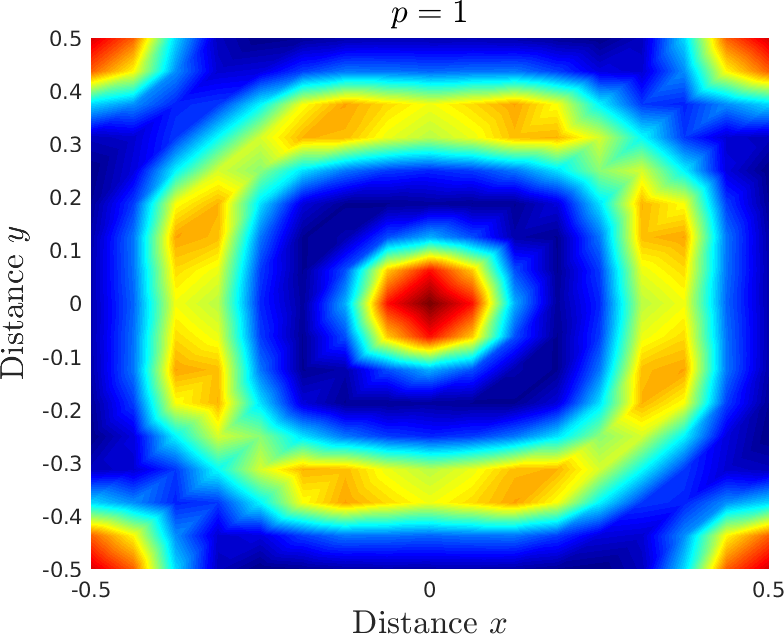}
\includegraphics[width=0.23\linewidth]{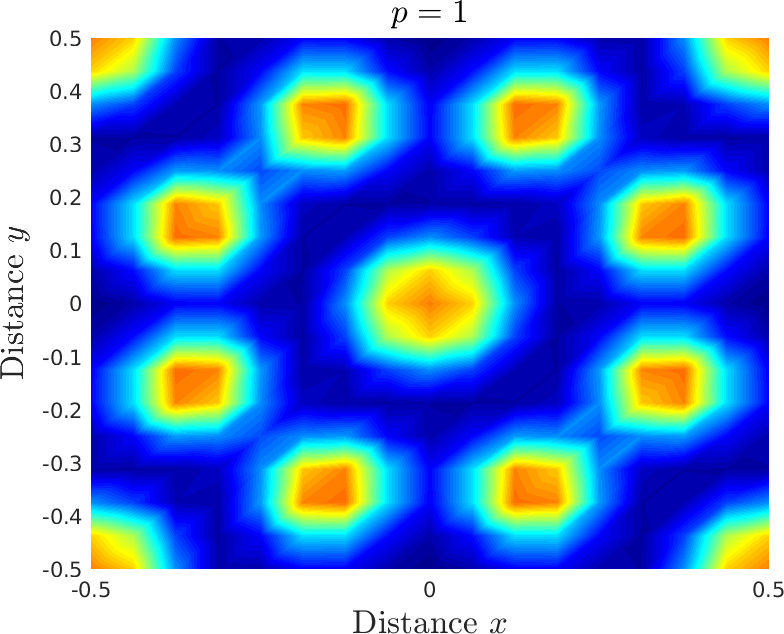}
\includegraphics[width=0.23\linewidth]{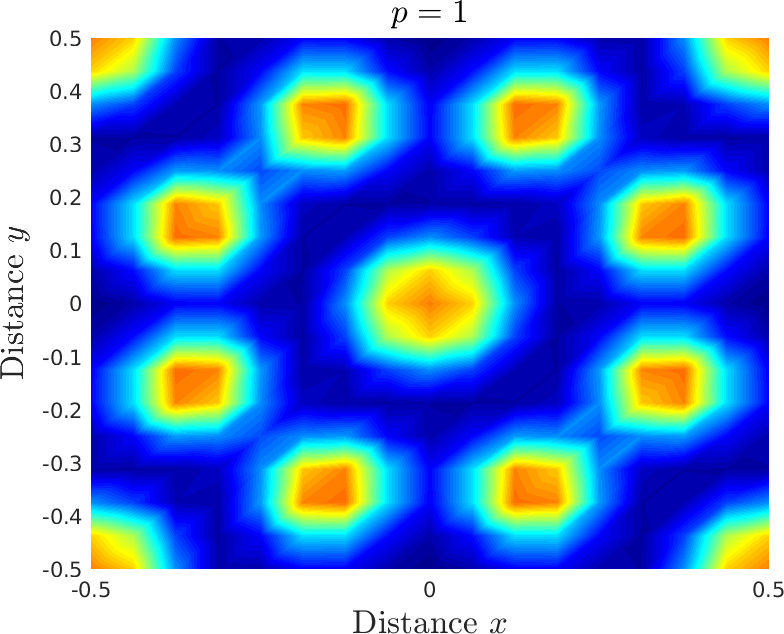}\\\vspace{0.1cm}
\includegraphics[width=0.23\linewidth]{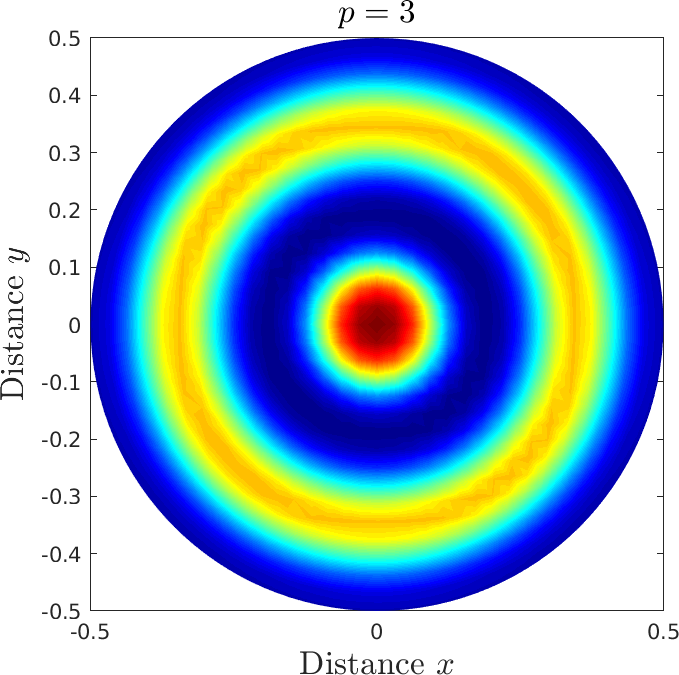}
\includegraphics[width=0.23\linewidth]{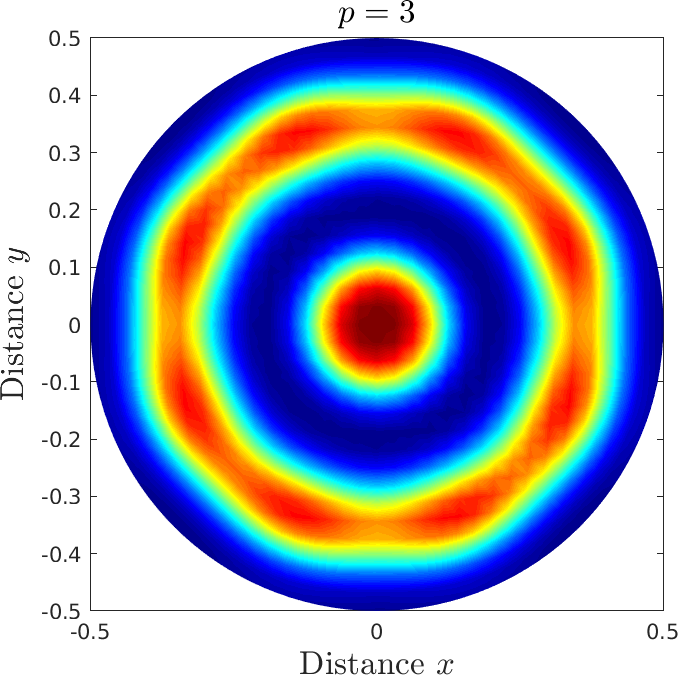}
\includegraphics[width=0.23\linewidth]{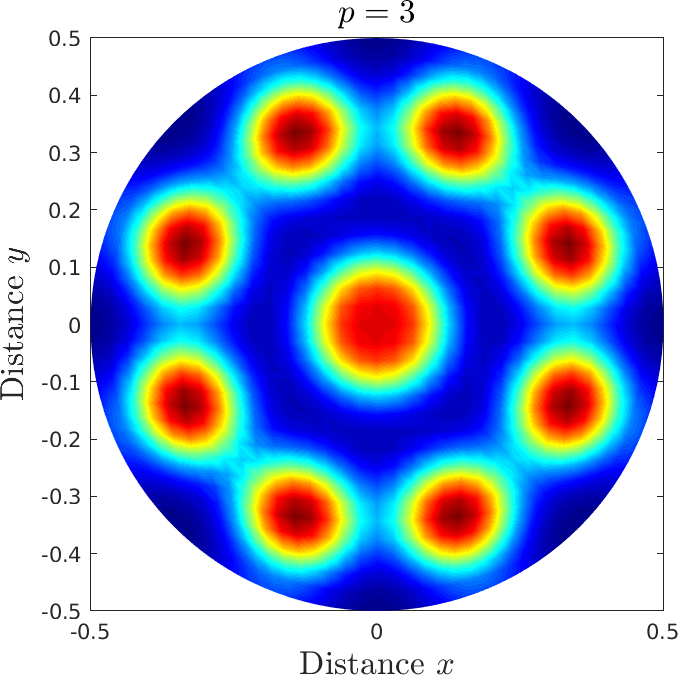}\\\vspace{0.1cm}
\includegraphics[width=0.23\linewidth]{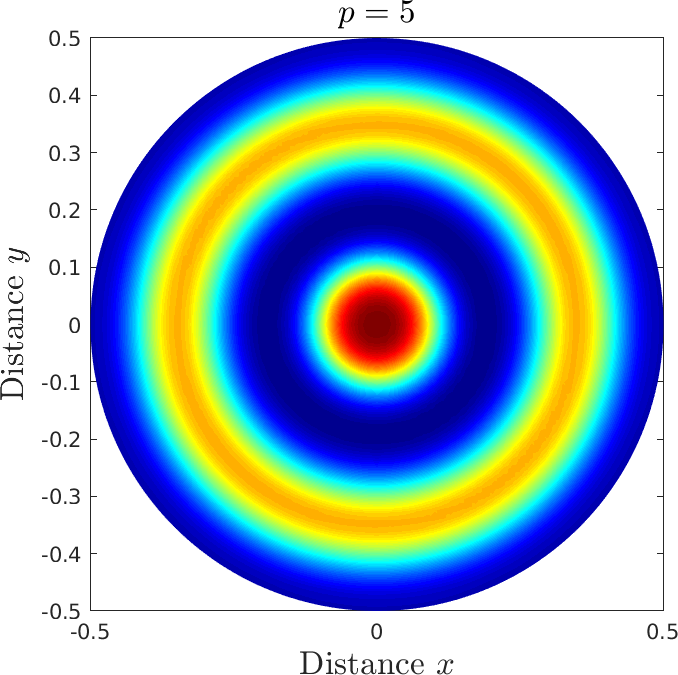}\includegraphics[width=0.23\linewidth]{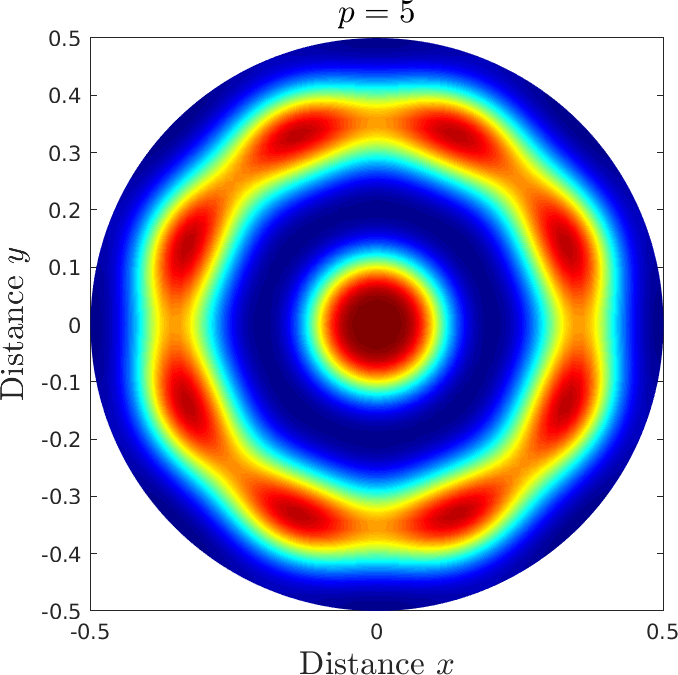}
\includegraphics[width=0.23\linewidth]{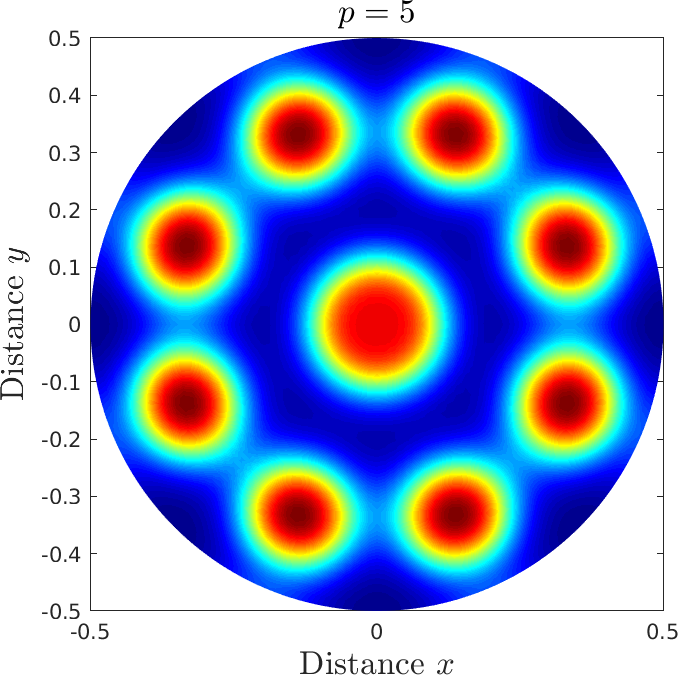}
\caption{Time evolution of the Turing patterns for different NURBS degrees at times $t = 0.5$, $t = 3$ and $t = 4$ using $\gamma = 100$.}\label{Ex2_Fig1}
\end{center}
\end{figure}
\begin{figure}[H]
\begin{center}
\includegraphics[width=0.23\linewidth]{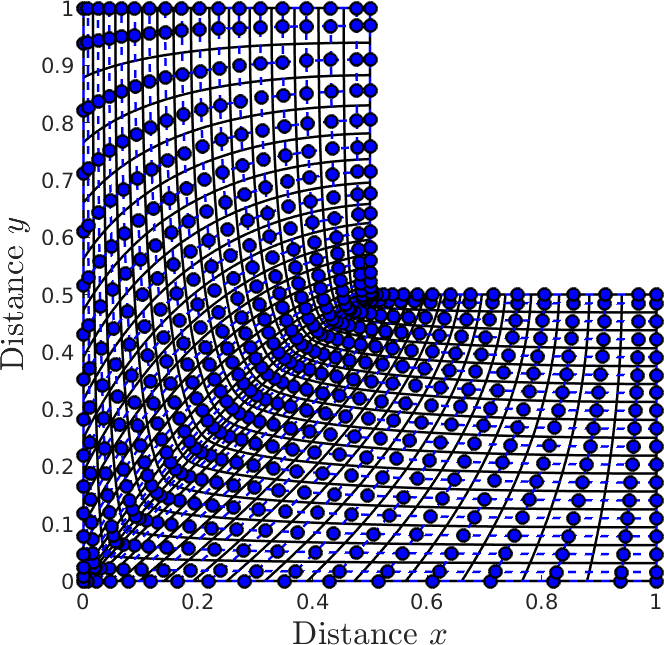}
\includegraphics[width=0.23\linewidth]{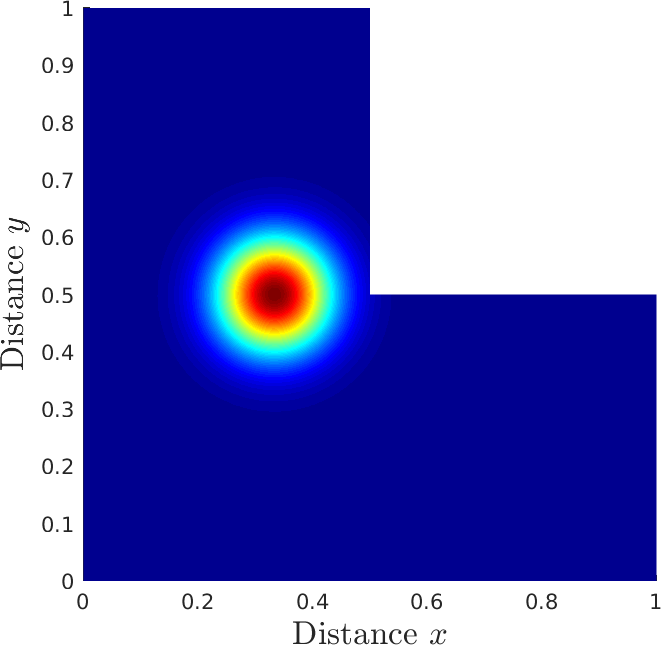}
\includegraphics[width=0.23\linewidth]{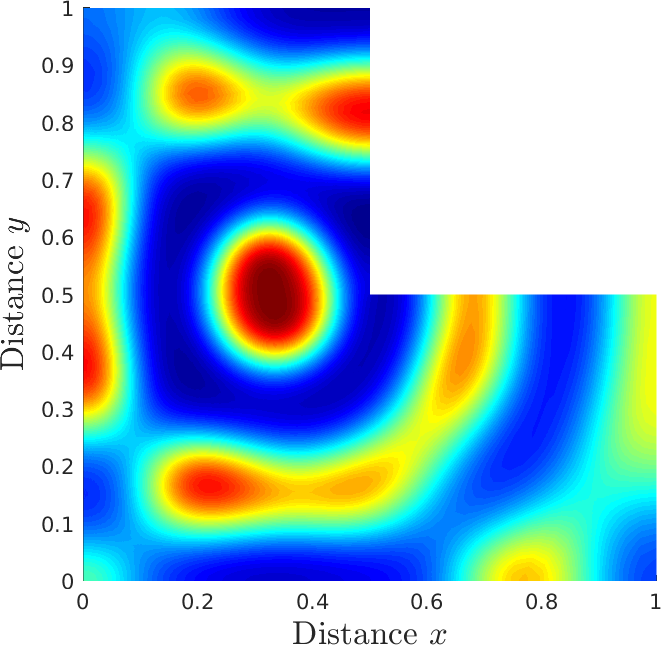}
\includegraphics[width=0.23\linewidth]{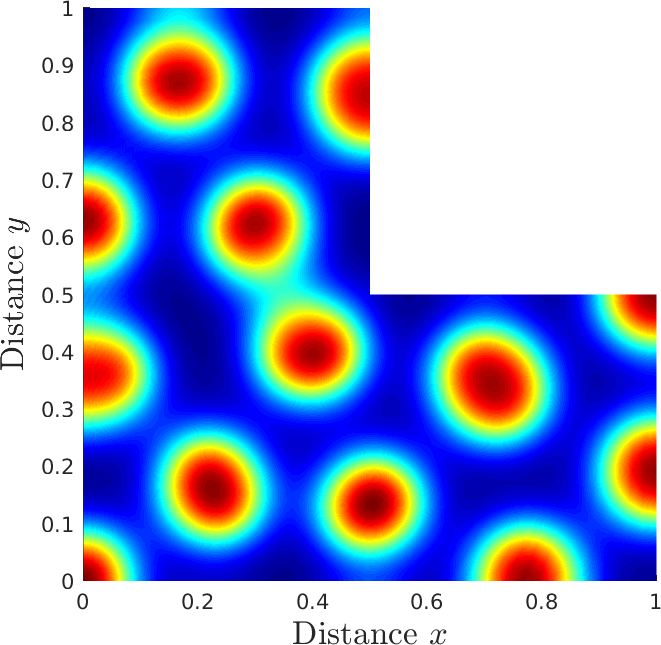}\\\vspace{0.1cm}
\includegraphics[width=0.23\linewidth]{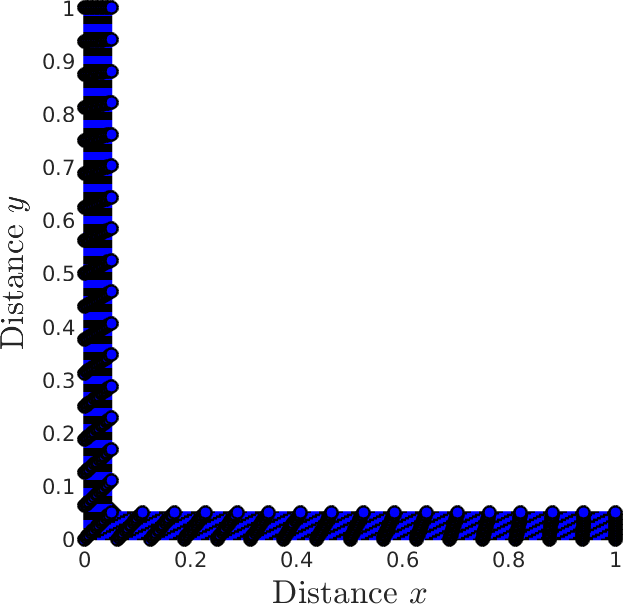}
\includegraphics[width=0.23\linewidth]{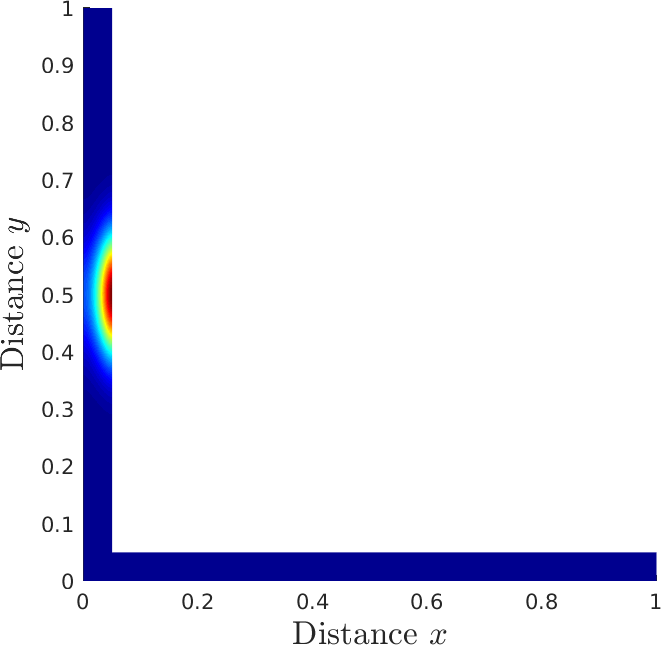}
\includegraphics[width=0.23\linewidth]{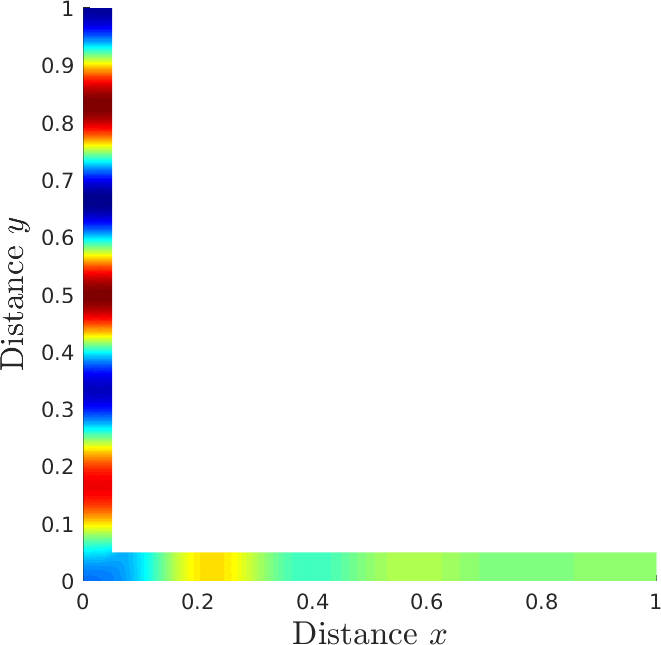}
\includegraphics[width=0.23\linewidth]{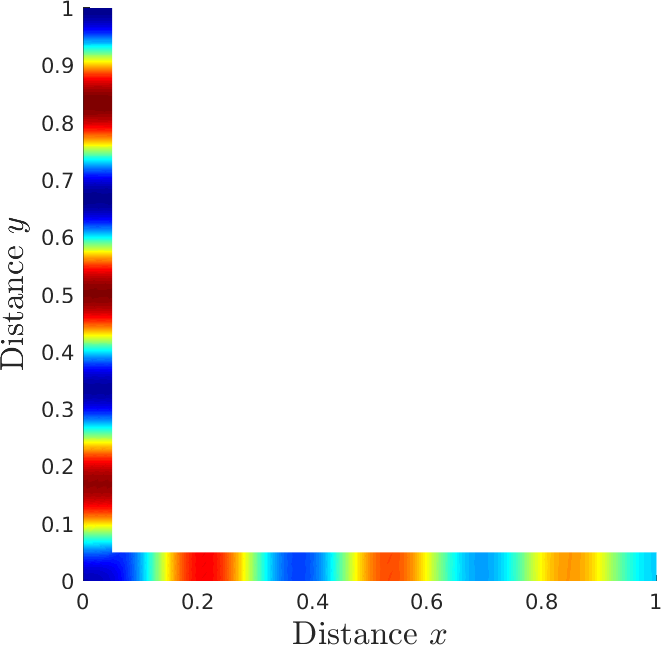}
\caption{Computational meshes and Turing patterns produced by the activator $u$ at times $t = 0, 0.5$ and $2$ on a wide L-shaped domain (top) and on a narrow L-shaped domain (bottom), using cubic NURBS.}\label{Ex2_Fig2}
\end{center}
\end{figure}
\subsection{Schnakenberg--Turing system}
\label{test3}
The Schnakenberg--Turing system \cite{turing1990chemical} is a nonlinear diffusion-reaction problem describing the growth of some embryos. The system consists of masses of tissue in which some substances react and diffuse in a process called morphogenesis. The Schnakenberg system is of the form \eqref{model}, where the reaction functions are given by
\[f(u,v) = \gamma(\alpha - u + u^{2}v),\quad g(u,v) = \gamma(\beta - u^{2}v),\]
where $\alpha$ and $\beta$ are positive constants representing the production of $u$ and $v$, respectively, while the nonlinear catalysis term $u^{2}v$ defines the activation of $u$ and the consumption of $v$. The initial conditions are chosen to be
\begin{eqnarray}
u(x,y,0) &=& \alpha + \beta + 10^{-3}e^{-100\left((x - \frac{1}{3})^{2} + (y - \frac{1}{3})^{2}\right)} ,\nonumber\\[-1ex]
\label{bur1}\\[-1ex]
v(x,y,0) &=& \frac{\beta}{(\alpha + \beta)^{2}}.\nonumber
\end{eqnarray}
Neumann boundary conditions are considered for this example. The parameter values are $\alpha = 0.1305$, $\beta = 0.7695$, while $\gamma$ is varied. The diffusion coefficients are $d_{1} = 0.05$ and $d_{2} = 1$. Our first objective in this test case is to investigate the performance of the novel IgA operator splitting method in the absence of advection, i.e., the pure diffusion-reaction system. Starting with linear NURBS, we discretize the computational domain given by the unit square $\Omega = [0,1]\times[0,1]$ into a uniform mesh and we perform a series of $k$-refinements to obtain high inter-element continuity. The resulting mesh resolution consists of a $32\times32$ square elements with element side length $h = 1/32$. Snapshots of the activator $u$ with $p = 1$, $p = 3$ and $p = 5$ are shown in Figure \ref{Ex2_Fig1}. Note that a rational quadratic basis is the minimum degree basis to represent a circular geometry. A first impression from Figure \ref{Ex2_Fig1} shows less accuracy achieved by the linear solution. We observe spurious oscillations over time. The solution is improved by using cubic NURBS, even though few oscillations are still detected (compare the cubic solution at time $t = 3$). By using NURBS of degree $p = 5$, high accuracy is achieved. The  initial spot is split over time into well-captured spots free of any oscillations and the computed results are in good agreement with those reported in the literature.\\

\begin{figure}
\begin{center}
\includegraphics[width=0.23\linewidth]{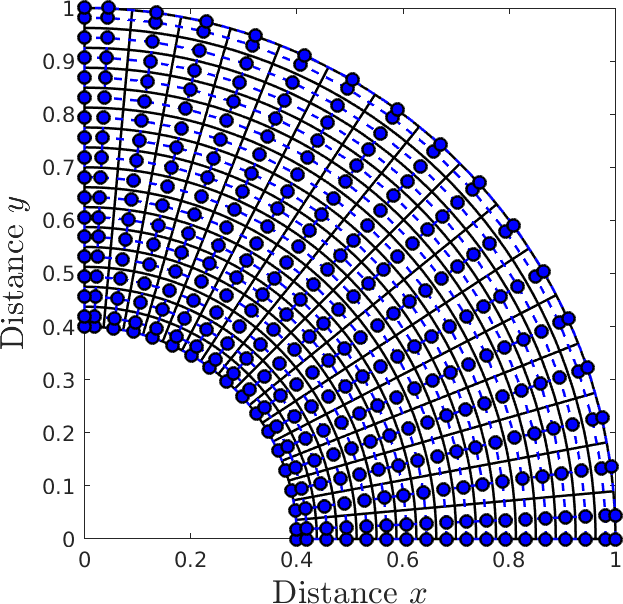}
\includegraphics[width=0.23\linewidth]{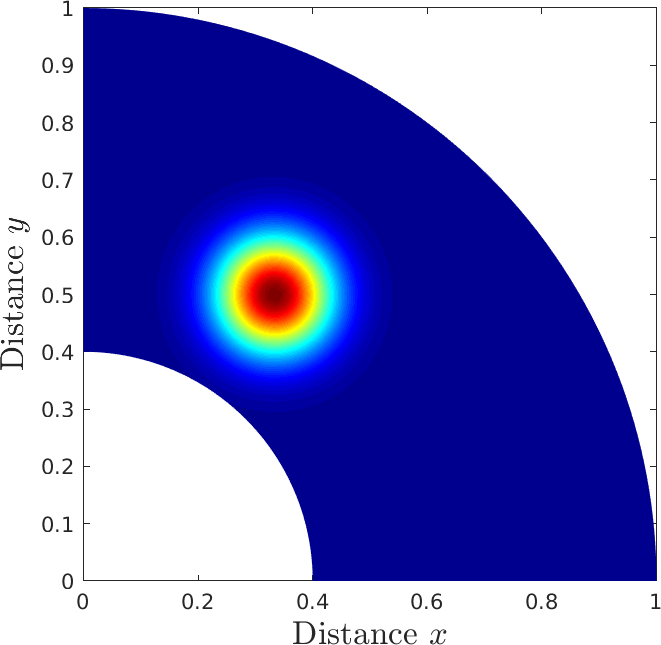}
\includegraphics[width=0.23\linewidth]{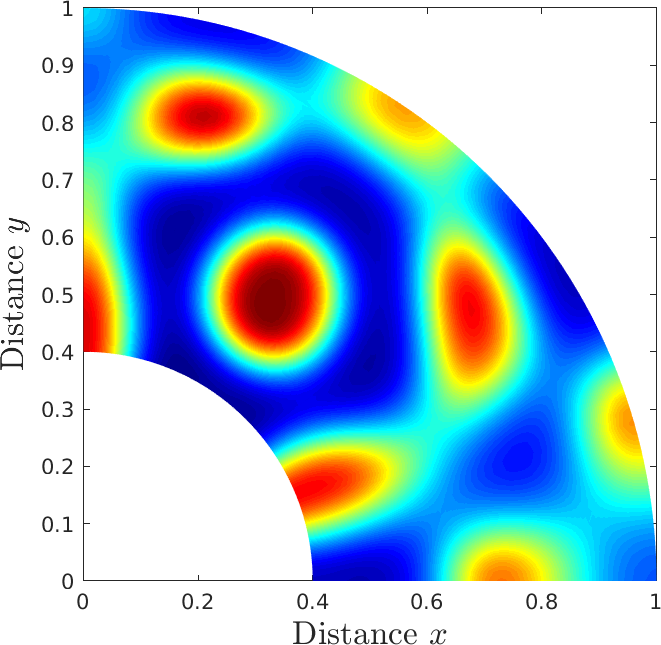}
\includegraphics[width=0.23\linewidth]{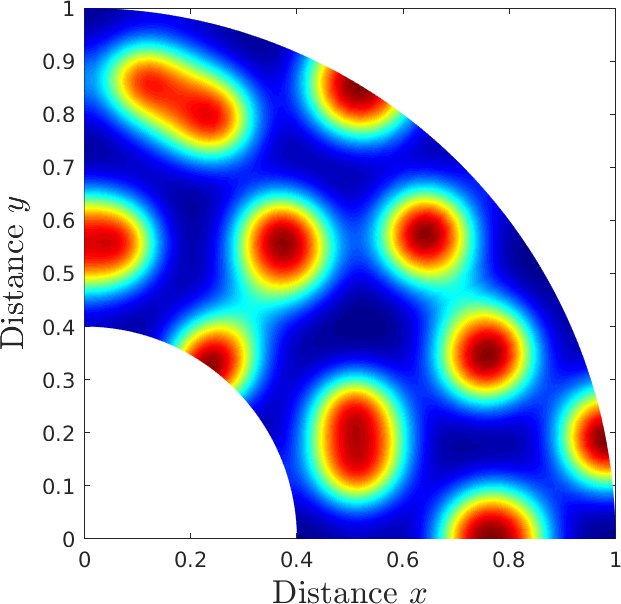}\\\vspace{0.1cm}
\includegraphics[width=0.23\linewidth]{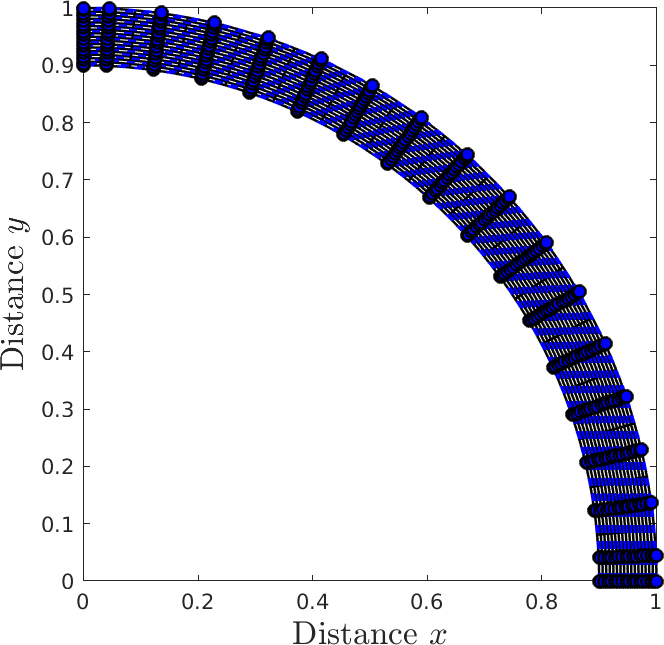}
\includegraphics[width=0.23\linewidth]{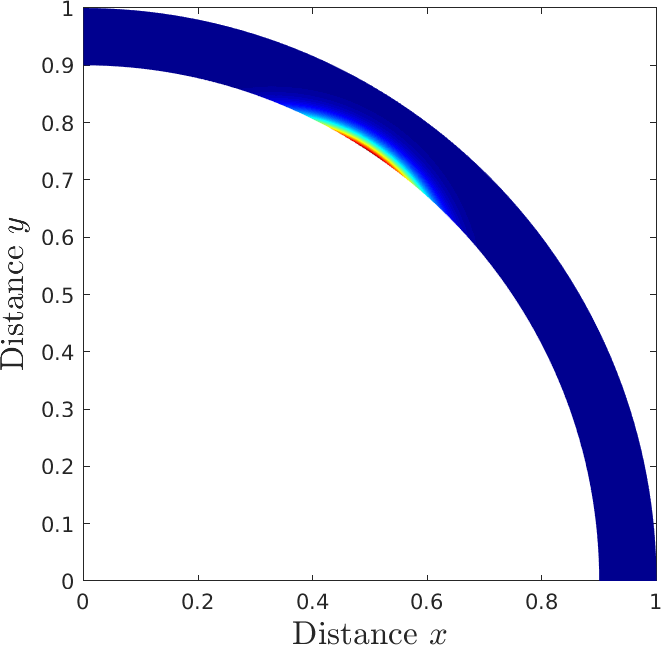}
\includegraphics[width=0.23\linewidth]{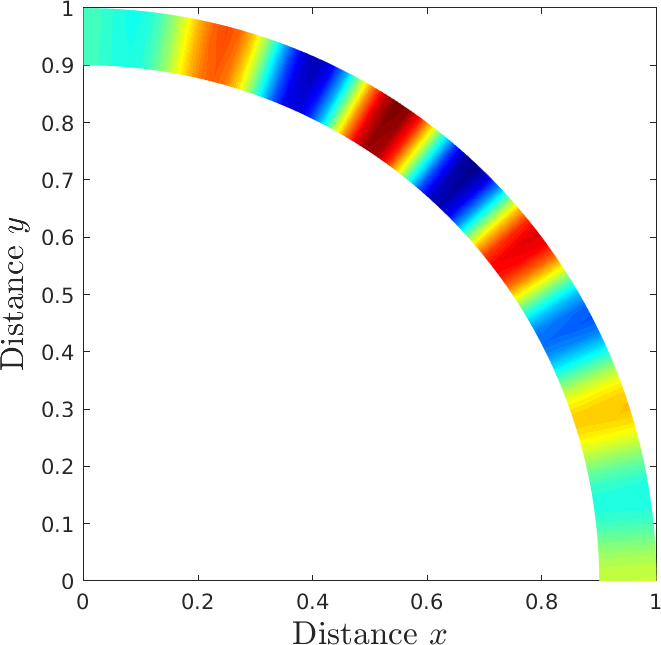}
\includegraphics[width=0.23\linewidth]{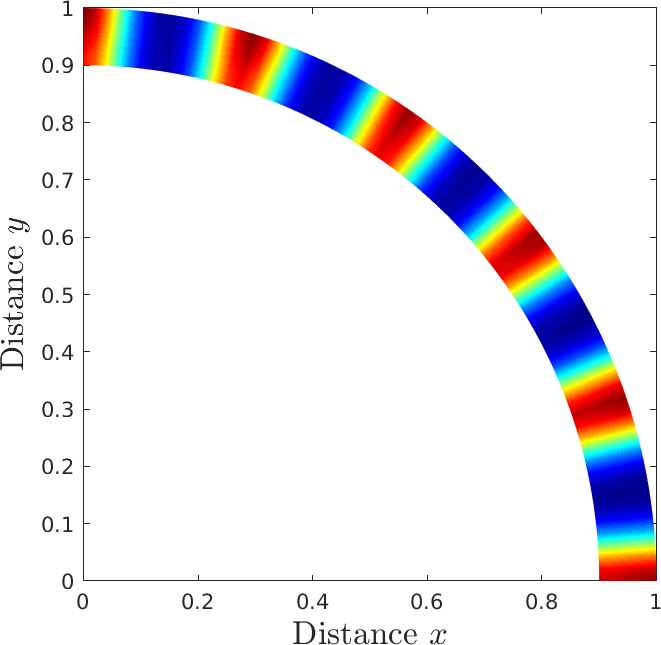}
\caption{Computational meshes and Turing pattern produced by the activator $u$ at times $t = 0, 0.5$ and $2$ on a wide annular domain (top) and on a narrow annular domain (bottom), using cubic NURBS.}\label{Ex2_Fig3}
\end{center}
\end{figure}

Next, to investigate the sensitivity of the Schankenberg--Turing problem with respect to complex geometries, two different geometries given by an $L$-shaped domain and an annular domain are explored. The initial condition again consists of a Gaussian that spreads out over time and replicates to form spot-like patterns. We keep all parameters at their previous values and fix $\gamma = 100$. We vary the inner radius of the annular domain and the diameter of the $L$-shaped domain to obtain narrow geometries as shown in the lower left of Figures \ref{Ex2_Fig2} and \ref{Ex2_Fig3}. Cubic NURBS are chosen for both geometries. We plot in Figures \ref{Ex2_Fig2} and \ref{Ex2_Fig3} the obtained snapshots of the activator $u$ using the $L$-shaped and the annular domain, respectively, at times $t = 0$, $t = 0.5$ and $t = 2$. It is interesting to see that $\gamma = 100$ results in spot-like patterns, in general. However, in the stiff narrow annular geometry and the narrow L-shaped domain, strip-like patterns appear; see the bottom row in Figures \ref{Ex2_Fig2} and \ref{Ex2_Fig3} respectively. This shows the sensitivity of Turing patterns to the geometry. A quantitative study confirms that our results are in good agreement with those published in \cite{zhu2009application}, where the same patterns are observed in a narrow rectangular geometry. This confirms the sensitivity of Turing patterns to geometry.\\


\begin{figure}
\begin{center}
\includegraphics[width=0.3\linewidth]{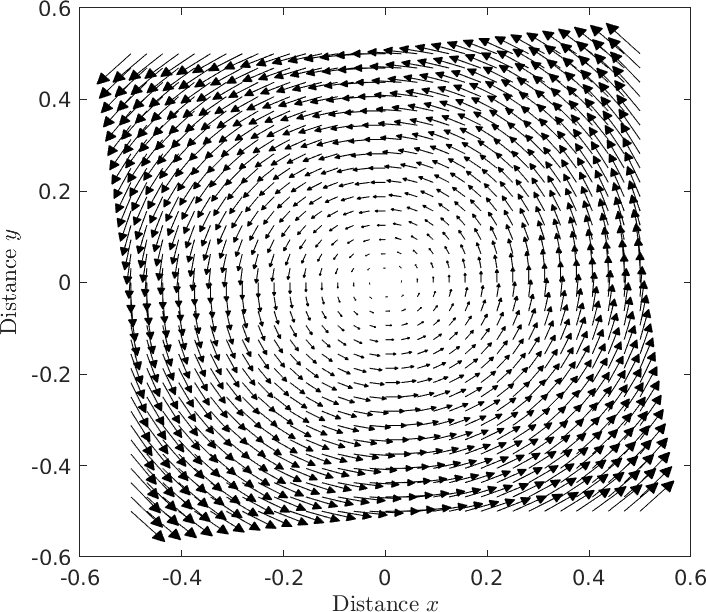}
\caption{Toroidal velocity field used for the computation of the Schnakenberg--Turing system.}\label{Ex2_Fig4}
\end{center}
\end{figure}
\begin{figure}
\begin{center}
\includegraphics[width=0.23\linewidth]{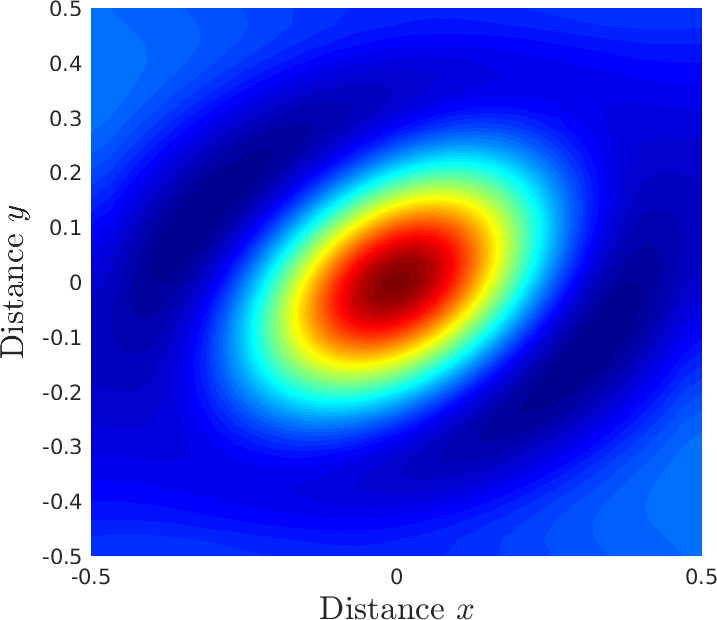}
\includegraphics[width=0.23\linewidth]{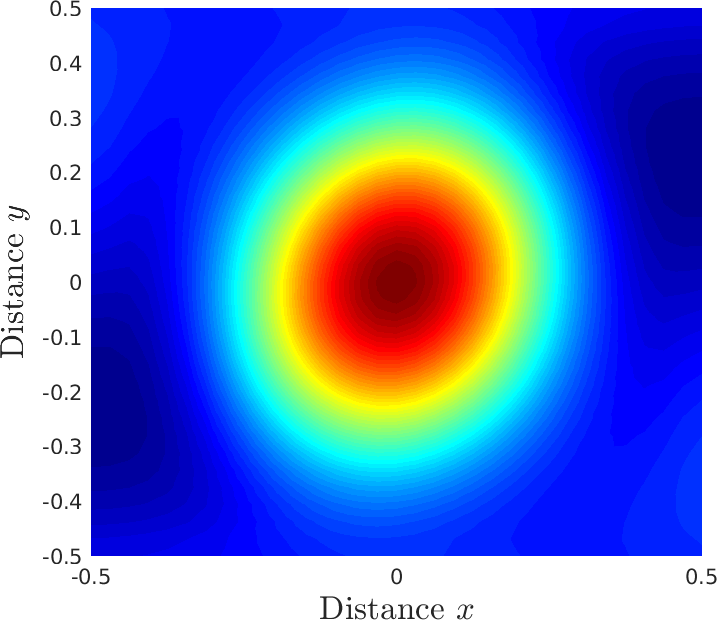}
\includegraphics[width=0.23\linewidth]{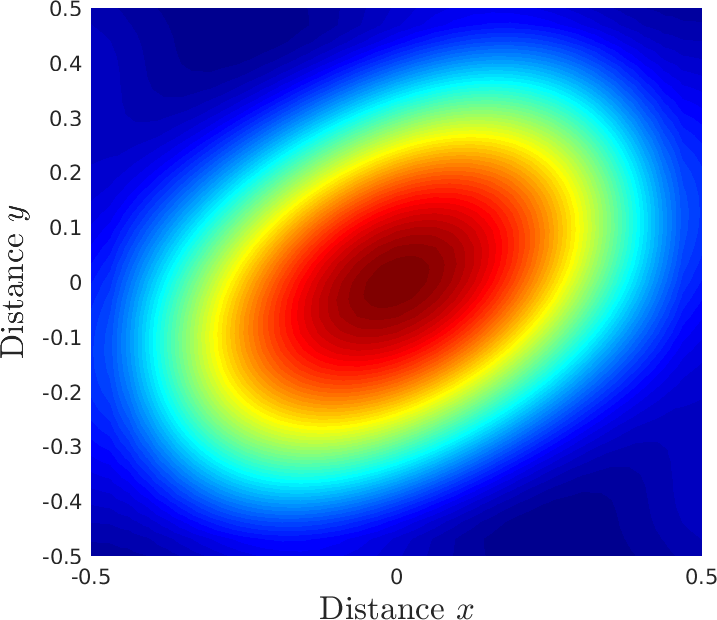}
\includegraphics[width=0.23\linewidth]{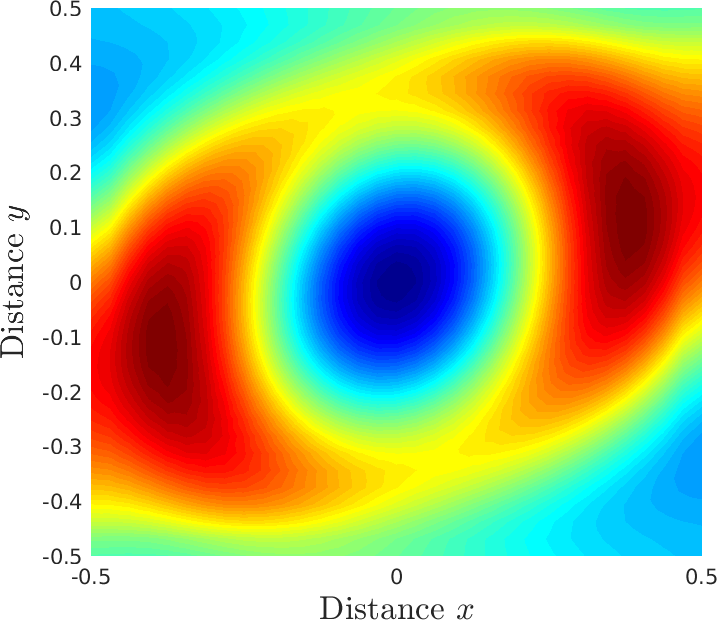}\\\vspace{0.1cm}
\includegraphics[width=0.23\linewidth]{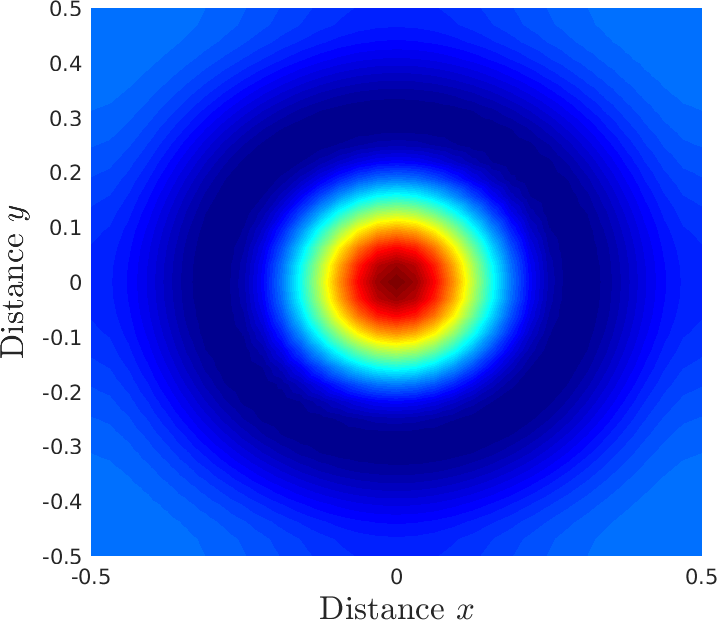}
\includegraphics[width=0.23\linewidth]{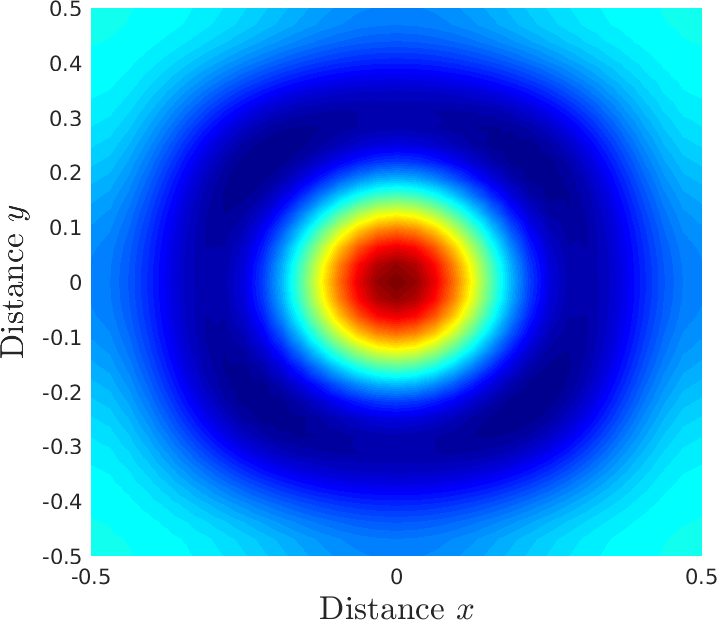}
\includegraphics[width=0.23\linewidth]{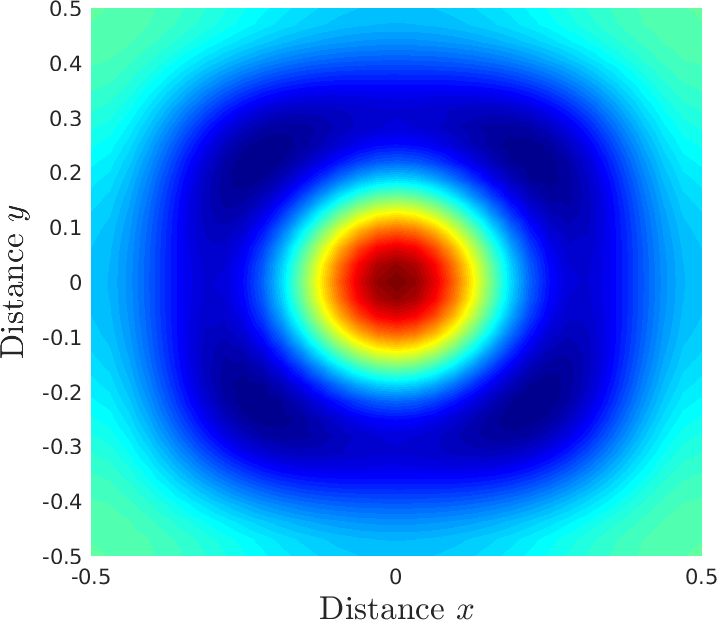}
\includegraphics[width=0.23\linewidth]{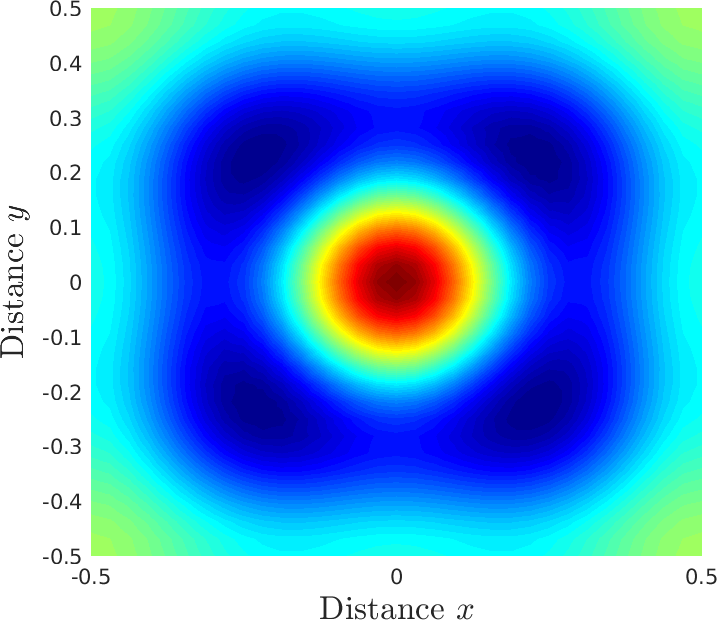}\\
\caption{Effect of the advection on the time evolution of the Turing patterns at times $t = 0.2$, $t = 0.4$, $t = 0.6$ and $t = 1$ using $\gamma = 30$. Top row: advection-diffusion-reaction case. Bottom row: diffusion-reaction case.}\label{Ex2_Fig5}
\end{center}
\end{figure}
\begin{figure}
\begin{center}
\includegraphics[width=0.23\linewidth]{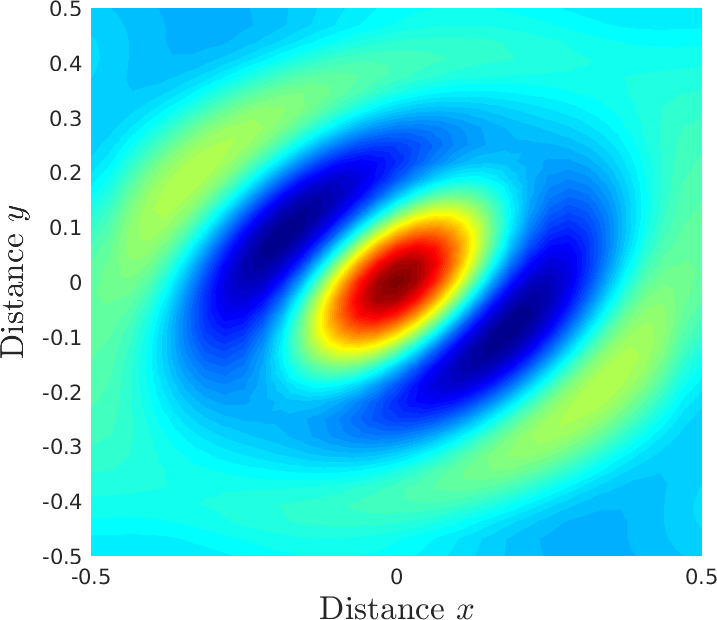}
\includegraphics[width=0.23\linewidth]{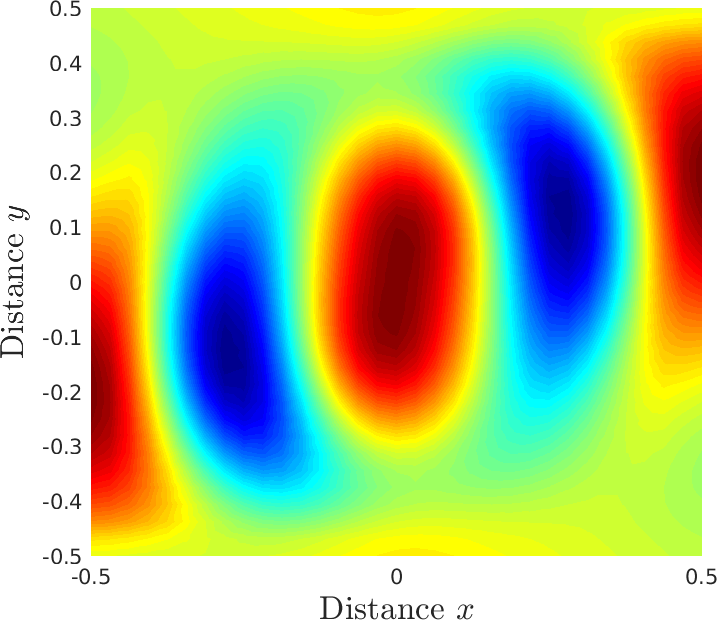}
\includegraphics[width=0.23\linewidth]{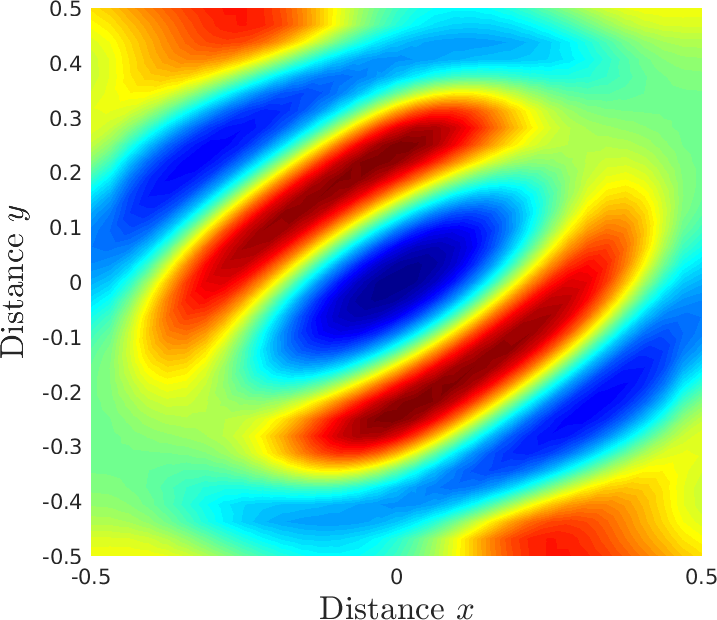}
\includegraphics[width=0.23\linewidth]{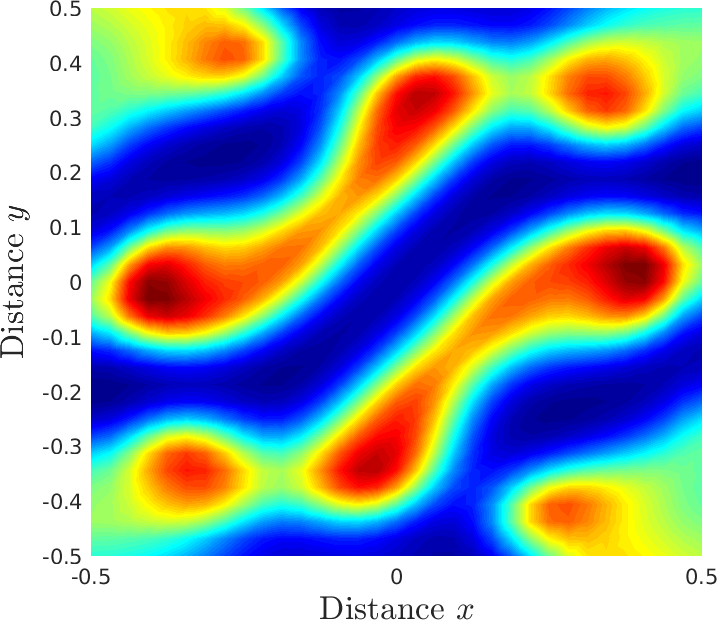}\\\vspace{0.1cm}
\includegraphics[width=0.23\linewidth]{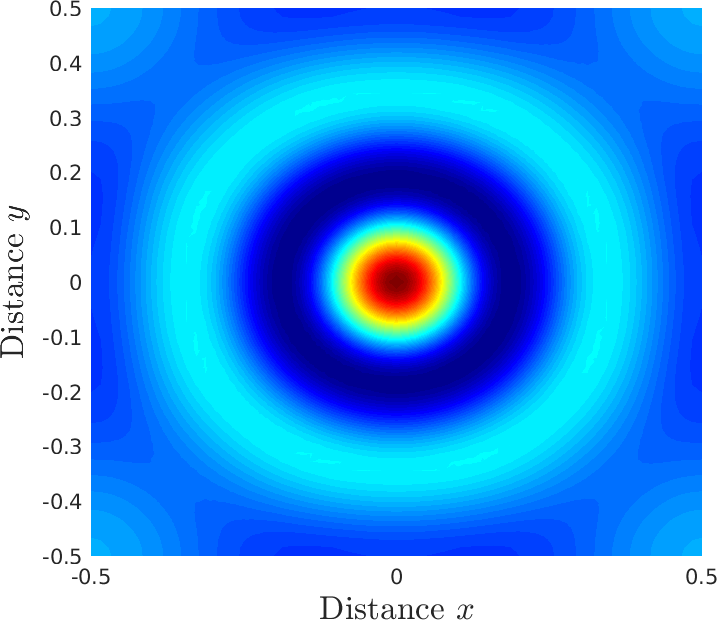}
\includegraphics[width=0.23\linewidth]{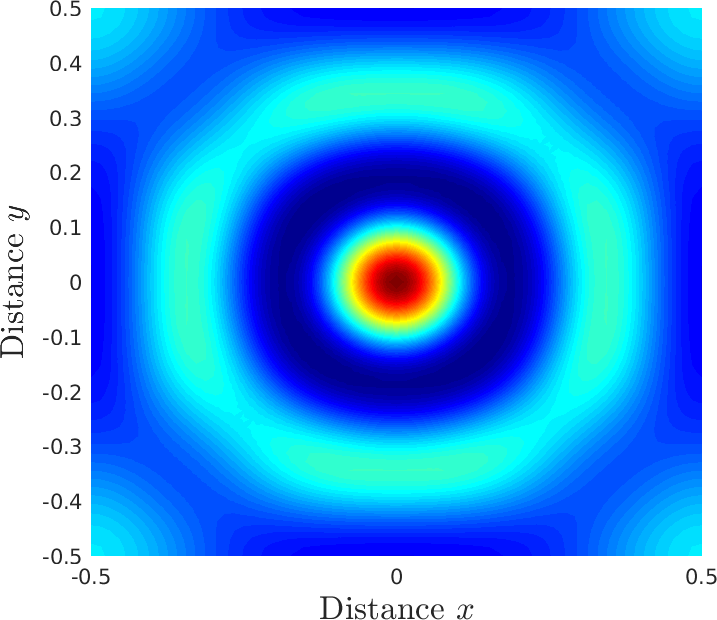}
\includegraphics[width=0.23\linewidth]{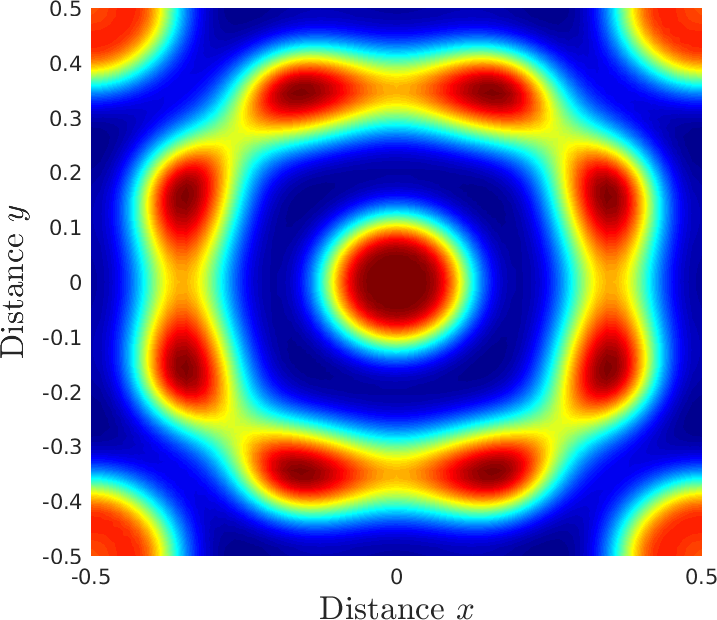}
\includegraphics[width=0.23\linewidth]{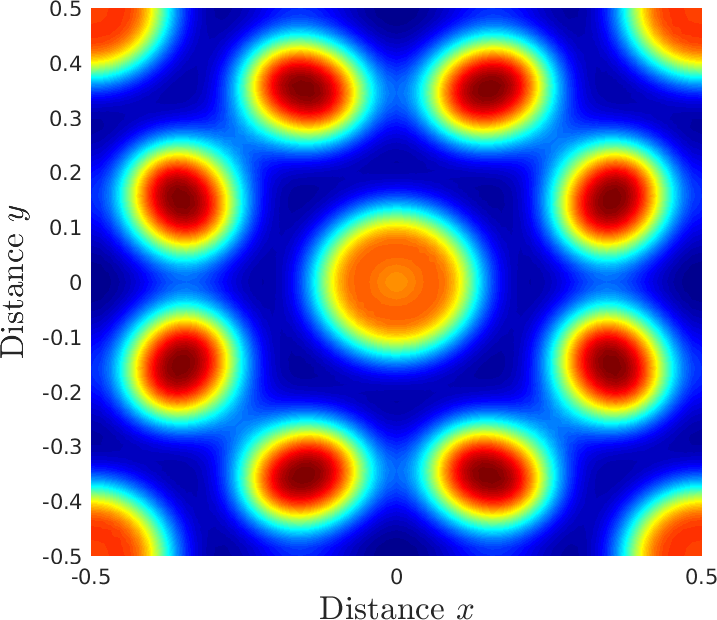}\\
\caption{Same as in Fig. \ref{Ex2_Fig5}, using $\gamma = 100$.}\label{Ex2_Fig6}
\end{center}
\end{figure}
\begin{figure}[!ht]
\begin{center}
\includegraphics[width=0.23\linewidth]{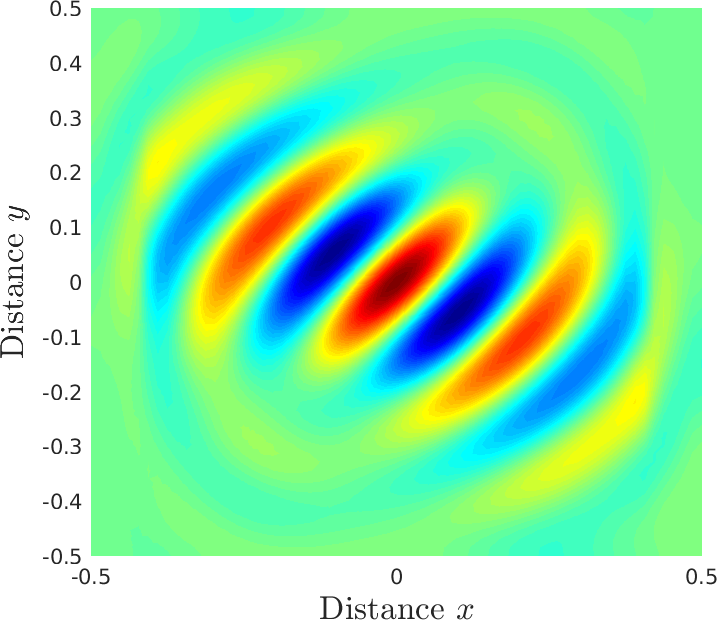}
\includegraphics[width=0.23\linewidth]{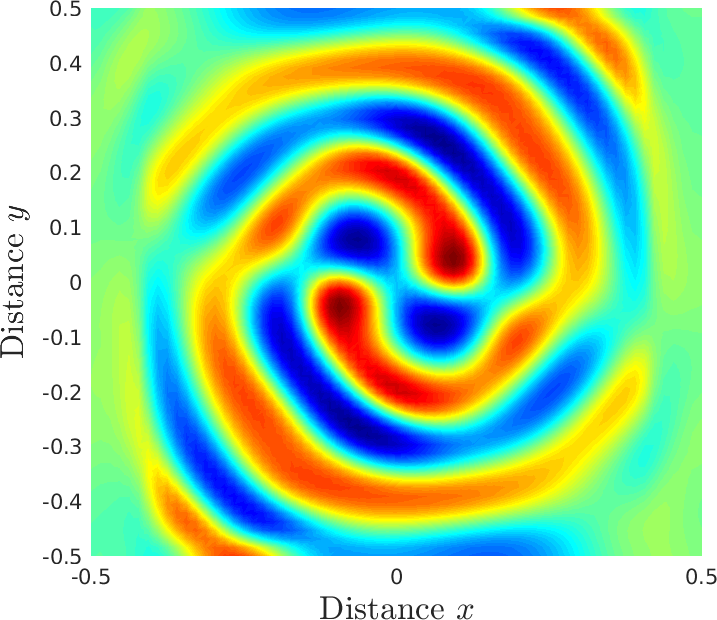}
\includegraphics[width=0.23\linewidth]{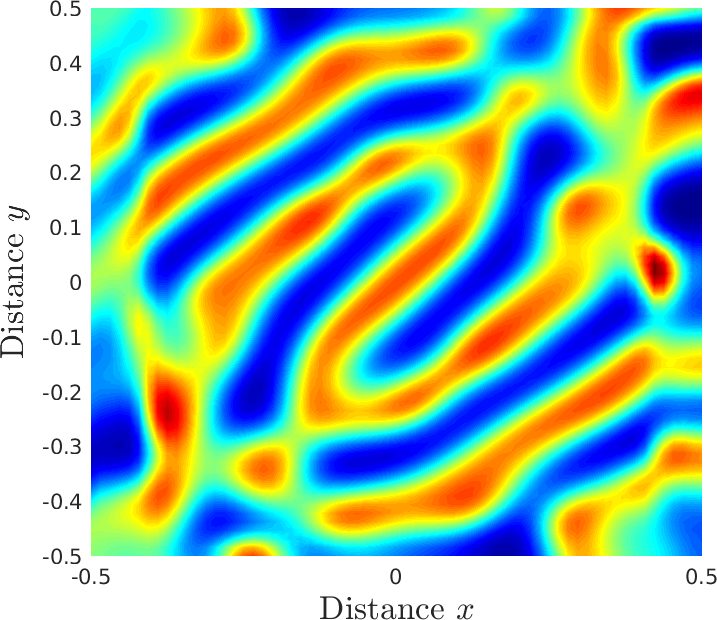}
\includegraphics[width=0.23\linewidth]{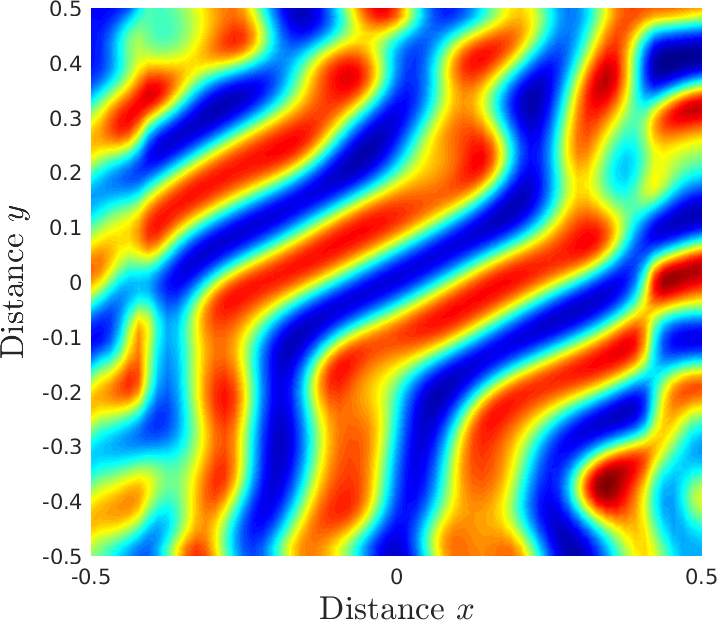}\\\vspace{0.1cm}
\includegraphics[width=0.23\linewidth]{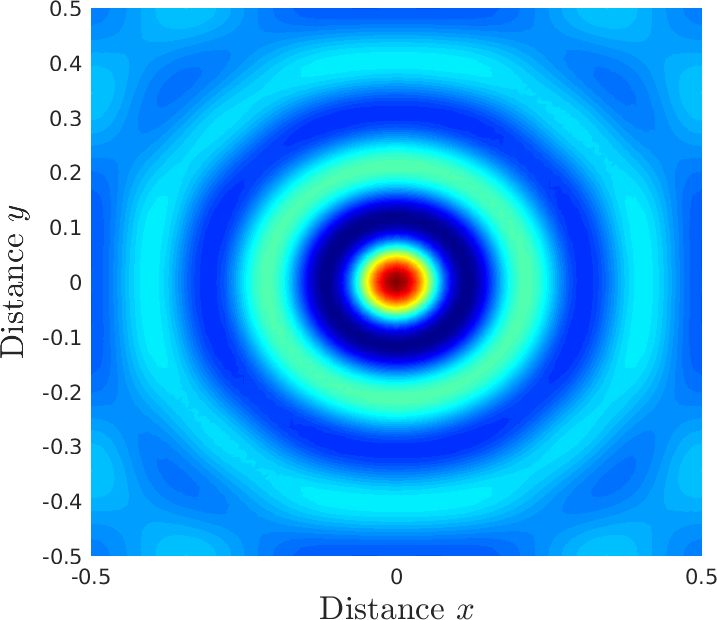}
\includegraphics[width=0.23\linewidth]{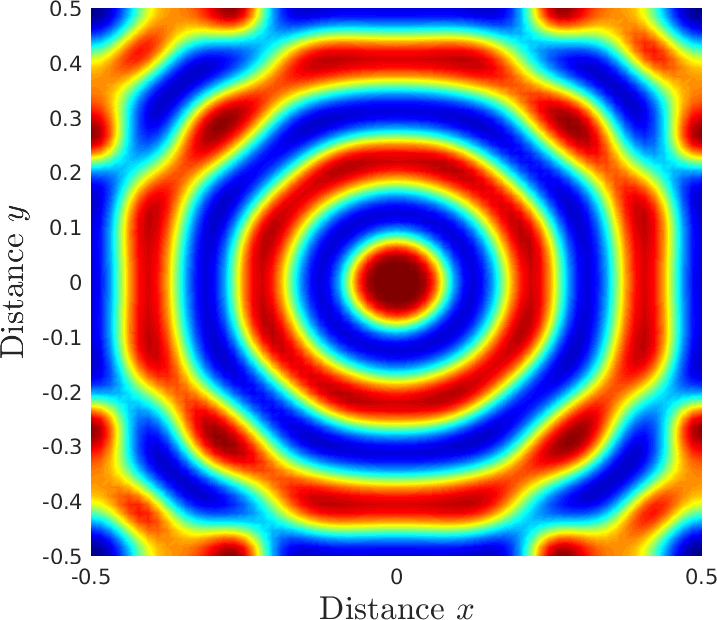}
\includegraphics[width=0.23\linewidth]{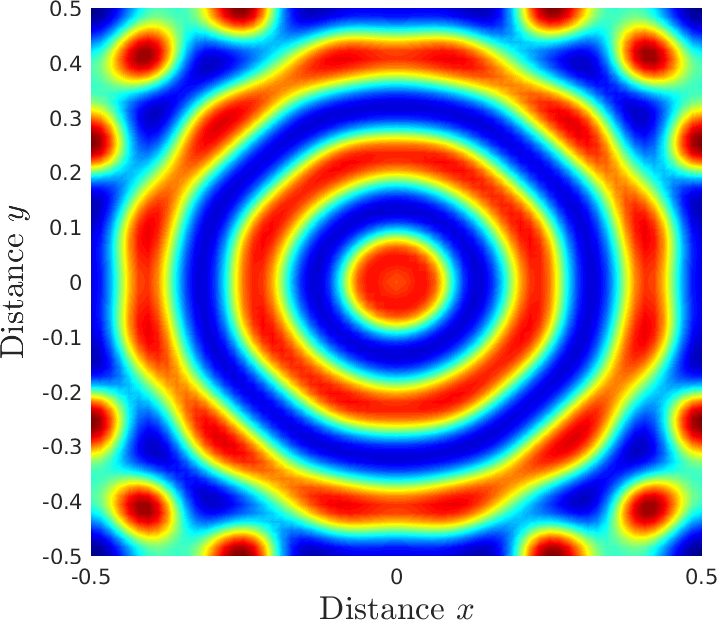}
\includegraphics[width=0.23\linewidth]{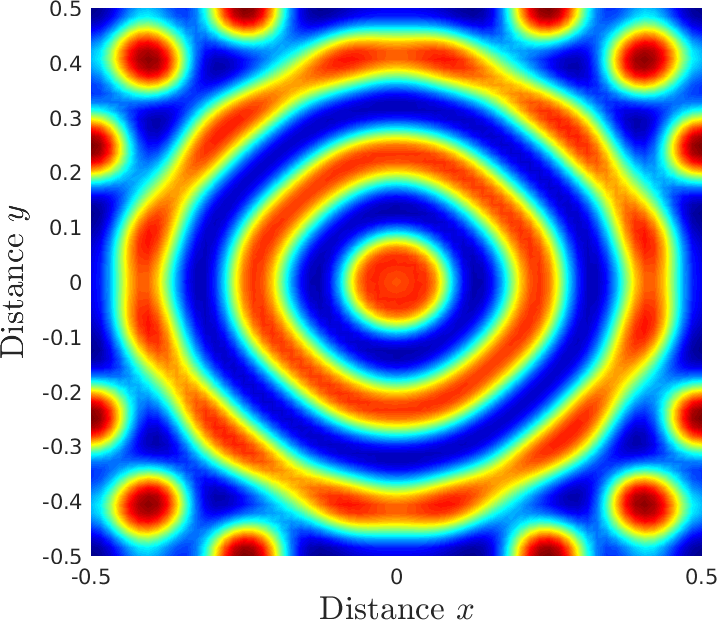}
\caption{Same as in Fig. \ref{Ex2_Fig5}, using $\gamma = 300$.}\label{Ex2_Fig7}
\end{center}
\end{figure}

Our next objective is to study the effect of the parameter $\gamma$ on Turing patterns in the presence of advection. Note that the parameter $\gamma$ measures the stiffness of the system. Indeed, the problem is considered stiff for values of $\gamma \geq 100$ \cite{zhu2009application,jiang2013krylov}. Moreover, the presence of advection adds an additional challenge to the system. Here, we choose a toroidal velocity field given by
\[{\bf a}({\bfx},t) = (-\omega y, \omega x),\]
where $\omega = 8$. The form of the velocity field is reported in Figure \ref{Ex2_Fig4}. The computational domain is given by the square $\Omega = [-0.5,0.5]\times[-0.5,0.5]$, which is divided into $32\times32$ elements using quartic NURBS. The behavior of the activator $u$ for various times and parameters $\gamma = 30$, $\gamma = 100$ and $\gamma = 300$ is given in Figures \ref{Ex2_Fig5},  \ref{Ex2_Fig6} and \ref{Ex2_Fig7}, respectively. The behavior of $u$ in the absence of advection is added to these figures to highlight the effect of the velocity field on Turing patterns. At time $t = 0.2$ with $\gamma = 30$, the initial spot located in the center of the domain is distorted by the direction of rotation of the toroidal velocity field, yielding two large spots at time $t = 1$. Note that in this case the reaction is negligible compared to diffusion, and diffusive effects are clearly observed. However, in the absence of advection, we observe that the initial spot given by the Gaussian remains almost the same during the time evolution, where the Gaussian shape is conserved. At the moderate value $\gamma = 100$, which can be considered as a critical value for which more spots are observed in the computed solution, the effects of the velocity lead to a deformation of the shape of the spots, which are distorted by the direction of rotation of the toroidal velocity field. It is interesting to observe that in the absence of advection, using $\gamma = 100$, a steady-state regime is almost reached at time $t= 1$, forming spot-like patterns; see the second row of Figure \ref{Ex2_Fig6}. However, the presence of advection continuously distorts the spots into arbitrary patterns; see the first row of Figure \ref{Ex2_Fig6}. At the high parameter value $\gamma = 300$, corresponding to the dominance of the reaction terms and resulting in strong stiffness, the formed strip-like patterns are again distorted by the direction of the toroidal velocity to turn around in the computational domain, as shown in Figure \ref{Ex2_Fig7}. On the other hand, the solution consists of distributed circular lines centered at the origin in the absence of advection. The addition of transport terms naturally distorts the original Turing patterns to produce patterns similar to the original, but with a certain deviation angle according to the magnitude of the velocity field.\\

\subsection{Gray--Scott system}\label{test3a}
\begin{figure}
\begin{center}
\includegraphics[width=0.3\linewidth]{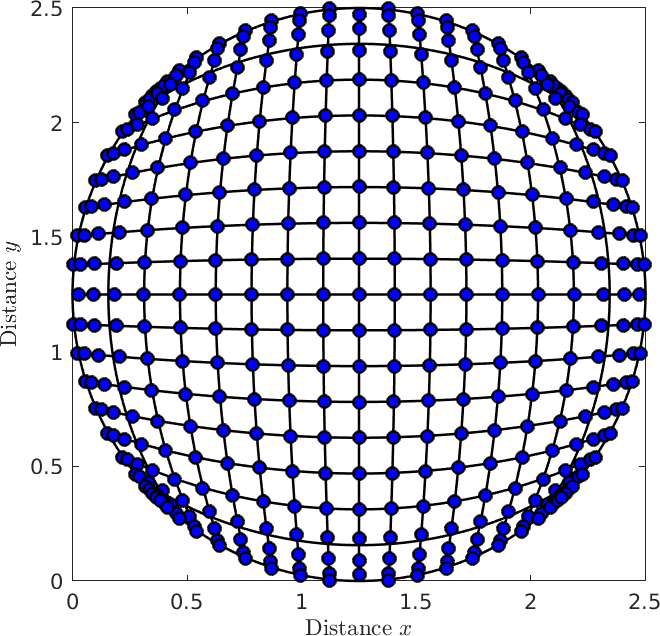}\hspace{2ex}
\caption{Computational domain used in the Gray--Scott system. Here, quintic NURBS in both $x$ and $y$ directions are used with 16$\times$16 elements and 441 control points.}\label{Ex3_Fig1}
\end{center}
\end{figure}

\begin{figure}[h!]
\begin{center}
\includegraphics[width=0.23\linewidth]{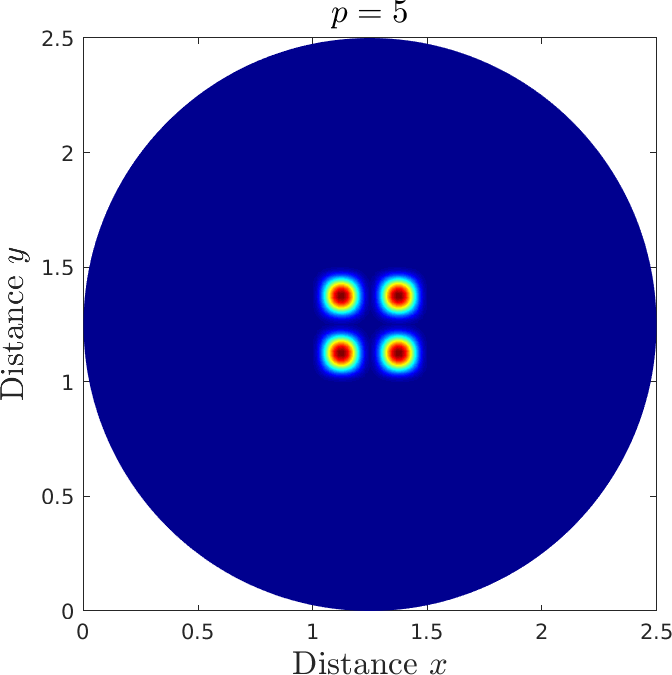}
\includegraphics[width=0.23\linewidth]{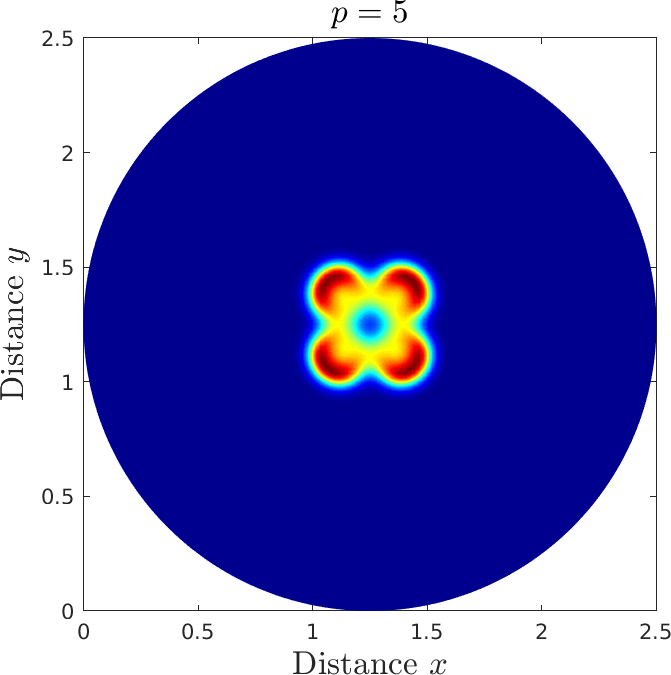}
\includegraphics[width=0.23\linewidth]{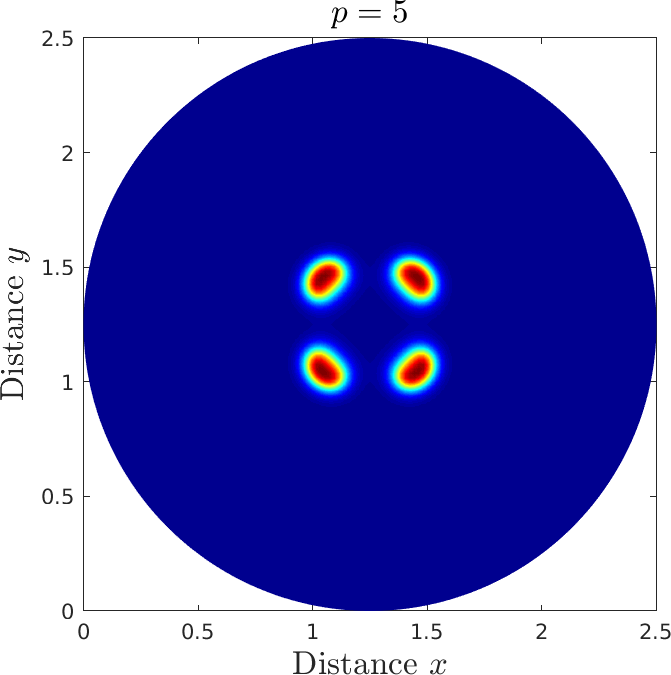}\\\vspace{0.1cm}
\includegraphics[width=0.23\linewidth]{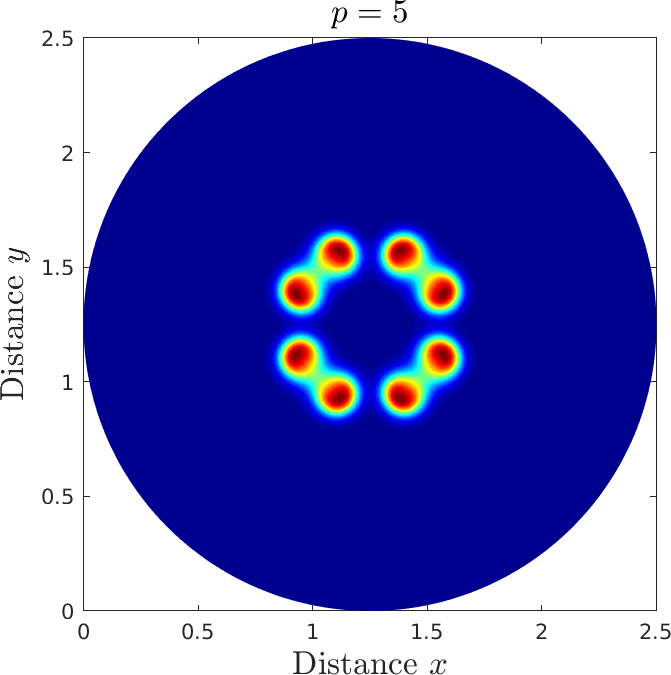}
\includegraphics[width=0.23\linewidth]{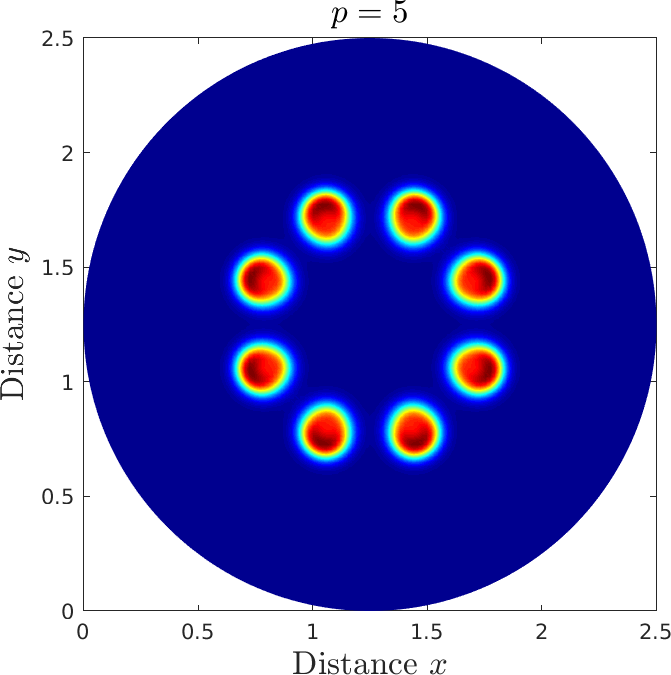}
\includegraphics[width=0.23\linewidth]{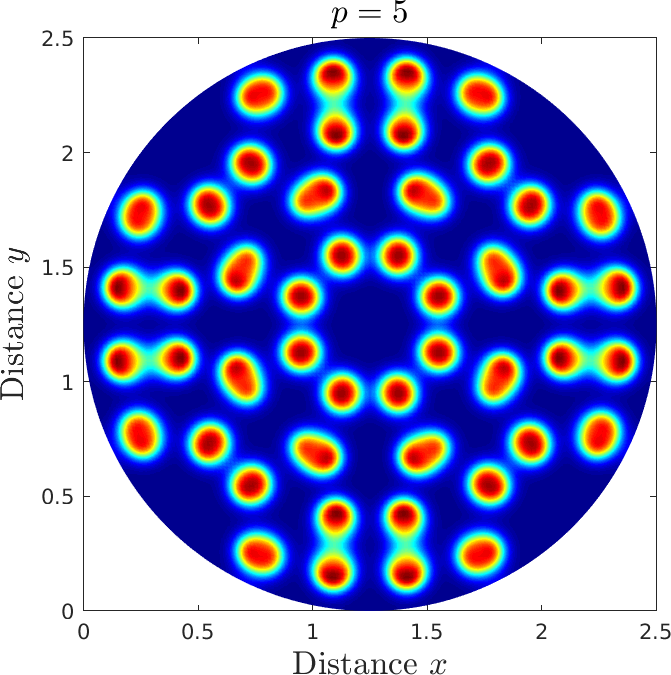}
\caption{Time evolution of the $v$ component in the Gray--Scott system at times $t = 1,100,200,500,1000,3000$, using quintic NURBS with 16$\times$16 elements and 441 nodes.}\label{Ex3_Fig2}
\end{center}
\end{figure}

The Gray--Scott system is given by problem \eqref{model}, where the reaction functions are given by
\[ f(u,v) = -uv^{2} + \alpha(1 - u),\quad g(u,v) = uv^{2} - (\alpha + \beta)v.\]
The parameters $\alpha$ and $\beta$ can be seen as bifurcation parameters \cite{delgado2017global}. Here, they are chosen as $\alpha = 0.024$, $\beta = 0.06$. Note that the parameter $\alpha$ refers to the rate at which the substance $u$ is replenished and the parameter $\beta$ controls the rate at which the substance $v$ is removed from the system during the reaction process. The diffusion coefficients are $d_{1} = 8.10^{-5}$, $d_{2} = 4.10^{-5}$. Typically, Neumann or periodic boundary conditions are used for this model. We have chosen Neumann conditions in the current study. The initial conditions are given by
\begin{eqnarray}
v(x,y,0) &=& \left\{
    \begin{array}{ll}
        \frac{1}{4}\sin^{2}(4\pi x)\sin^{2}(4\pi y) & \mbox{if } \{x,y\} \in [1,1.5] \\
        0 & \mbox{else}
    \end{array}
\right. ,\nonumber\\[-1ex]
\label{bur1a}\\[-1ex]
u(x,y,0) &=& 1 - 2v(x,y,0).\nonumber
\end{eqnarray}
Using these initial conditions, the solution consists of four spots centered in the domain. During the time evolution, the initial four spots will split and multiply into spot-like patterns. This problem has previously been analyzed using a moving finite element method \cite{hu2012moving} and a moving element free approach \cite{dehghan2016numerical} on a square domain given by $[0,2.5]\times[0,2.5]$. In this work, the Gray--Scott model is solved on a circular domain centered at $(x_{0},y_{0}) = (1.25,1.25)$. The spatial mesh resolution is set to $16\times16$ elements using quintic NURBS functions as shown in Figure \ref{Ex3_Fig1}. In Figure \ref{Ex3_Fig2}, we present the snapshots of the $v$-component at different times. Clearly, the spots are well captured and no oscillations are observed. In comparison with the results published in \cite{dehghan2016numerical,hu2012moving}, we have reproduced the same patterns and our results are in good agreement with those in \cite{dehghan2016numerical,hu2012moving}. It is worth noting that the novel IgA operator splitting method successfully captures the patterns reproduced in the Gray--Scott system using a very coarse NURBS mesh.
\section*{Conclusion}
In this work we have presented a new IgA operator splitting method for solving a class of coupled advection-diffusion-reaction equations. The new approach is based on separating the stiff part of the problem from the non-stiff one and solving each part with an appropriate technique. The advective terms, which are usually the source of numerical instabilities, are treated in a semi-Lagrangian framework where the NURBS functions together with an $L^{2}$-projection approach are used to update the solution during the advection stage. The resulting diffusion-reaction equations, described in terms of material derivatives, are then split into a linear diffusion part and a nonlinear reaction part, based on the Strang splitting technique to obtain second-order accuracy. The diffusive terms are discretized using IgA, and the resulting ODEs are solved using an implicit second-order scheme, while the reactive terms are integrated using an explicit scheme. The main objective of this study was to reduce the complexity resulting from the multiphysical nature of the overall system, which was successfully achieved by using Strang operator splitting in combination with well-established algorithms to solve the subproblems. Although the new method has many similarities with standard finite element analysis, the $k$-refinement strategy implemented in the present study has no analogy in conventional finite element methods. The key feature of this technique is essentially to ensure high inter-element continuity by providing smoother functions than those of conventional finite element methods. This is due to the fact that the order elevation process ($p$-refinement) and the knot insertion process ($h$-refinement) do not commute, see \cite{hughes2005isogeometric,da2011some}. The computed results support our expectations for a stable and highly accurate new IgA operator splitting method for nonlinear advection-diffusion-reaction problems. It is worth noting that the $k$-refinement strategy adopted in this paper shows similar behavior for the selected test cases. The novel method also avoids linearization of the reaction terms. Keeping reaction and diffusion together can be computationally demanding. The presented numerical results are in good agreement with existing studies. Moreover, by the addition of transport terms, the original Turing patterns were distorted. Patterns similar to the original ones were created, but with a certain deviation angle which was proportional to the magnitude of the velocity field.
\bibliographystyle{elsarticle-num-names}
\bibliography{refs}
\end{document}